%% file: verifiable_conditions_manifolds.tex
\documentclass[12pt]{article}

\usepackage{authblk}
\usepackage[margin=2cm]{geometry}

\input{header}
\input{def_armand}

\input{def}

\usepackage{xr-hyper}

\begin{document}

		\title{On the irreducibility and convergence of a class of nonsmooth nonlinear state-space models on manifolds and their applications to zeroth-order optimization}  
		
		
		\author[1,2]{Armand Gissler} 
		\author[2]{Alain Durmus}
		\author[1,2]{Anne Auger}
		
		\affil[1]{Inria Saclay Île-de-France, Palaiseau, France}
		\affil[2]{CMAP, CNRS, École polytechnique, Institut Polytechnique de Paris, Palaiseau, France}
		
		\maketitle
		
	\begin{abstract}
		
		In this paper, we analyze a large class of general nonlinear state-\break space models on a state-space $\cX$, defined by the recursion $\phi_{k+1}  = \break F(\phi_k,\alpha(\phi_k,U_{k+1}))$, $k \in\bN$, where $F,\alpha$ are some functions and $\{U_{k+1}\}_{k\in\bN}$ is a sequence of \iid~random variables.
		More precisely, we
		extend conditions under which this class of Markov chains is irreducible, aperiodic and satisfies important continuity properties, 
		relaxing two key assumptions from prior works. 
		First, the state-space $\cX$ is supposed to be a smooth manifold instead of an open subset of a Euclidean space.
		Second, we only suppose that $F$ is locally Lipschitz continuous.

		We demonstrate the significance of our results through their application
		to Markov chains underlying optimization algorithms. 
		These schemes belong to the class of evolution strategies with covariance matrix adaptation and step-size adaptation.
	\end{abstract}
		\newcommand{\kwd}[1]{#1}
		
		

\section{Introduction}
Consider a nonlinear state-space model defined by the recursion:
\begin{equation}\label{eq:NLstate-space-model}
	\phi_{k+1} = G(\phi_k, \xi_{k+1})\eqsp,
\end{equation}
where the sequence $\{\xi_{k+1}\}_{k\in\mathbb{N}}$ consists of independent and identically distributed (i.i.d.) random variables,
$G: \msx \times \msw \to \msx$ is a continuous function, and $\msx,\msw$
are two measurable spaces.
Nonlinear state-space models \eqref{eq:NLstate-space-model} form
a class of Markov chains that have been first popularized in stochastic control theory 
\cite{meyn1993model,meyn1987new,meyn1989stochastic,meyn1991asymptotic}. This has spurred extensive analysis and has a well-established historical context. In particular, for nonlinear autoregressive models, i.e., where $G$ can be written as $G(x,u ) = \tilde{G}(x) + u$,
ergodicity has been widely investigated \cite{bhattacharya1995geometric,an1997note,yao2000stability,mokkadem1987criteres,glynn2017recurrence}.
Moreover, connections have been established between the stability of \eqref{eq:NLstate-space-model} and the one of some Ordinary Differential Equation (ODE) \del{in }\cite{huang2002ode}.
The idea of analyzing \eqref{eq:NLstate-space-model} from the perspective of control theory, where $u$ is regarded as a control parameter, was initially proposed in \cite{stroock1972support} within the context of diffusion processes. This approach was subsequently employed with success in \cite{ichihara1974classification} and \cite{kliemann1987recurrence}.
It has been then applied in \cite{meyn1989stochastic,meyn1991asymptotic} when $G$ is infinity differentiable and $\msx,\msw$ are open sets of Euclidean spaces, to establish the irreducibility, aperiodicity and topological properties of the Markov kernel associated to \eqref{eq:NLstate-space-model}. {The theory developed in \cite{meyn1989stochastic,meyn1991asymptotic} forms the foundation of \cite[Chapter~7]{meyn2012markov}, which in turn underpins the present work.}

Besides stochastic control models, \eqref{eq:NLstate-space-model} also encompasses many algorithms in
optimization and Markov chain Monte Carlo algorithms. The variable $\phi_\t$
corresponds to the state of an algorithm at iteration $\t$ and
$\del{U}{\xi}_{\t+1}$ represents the random components used to update this
state. 
However, for certain classes of algorithms, especially those arising from zeroth-order optimization, the function $G$ may not be continuous. 
Nevertheless, they can be written using the Markov chain model introduced in \cite{chotard2019verifiable} as 
\begin{equation}
	\label{eq:def-MC}
	\phi_{\t+1} = F(\phi_\t,\alpha(\phi_\t,U_{\t+1}))\eqsp,
\end{equation}
for some continuous function $F: \msx \times \msw \to \msx$ and a potentially discontinuous function
$\alpha : \msx \times \msu \to \msw$ assumed measurable,
and $\seq{U_{\t +1} }{\t\in\bN}$ a sequence of
\iid~random variables valued in a measurable space
$(\msu,\mcu)$, chosen independently of the initial state $\phi_0$.
{When taking $\alpha(x,u)=u$, the model in \eqref{eq:def-MC} is equivalent to \eqref{eq:NLstate-space-model}. Leveraging the form \eqref{eq:def-MC}, under the assumption that $F$ is continuously differentiable and under suitable assumptions on $\alpha$ (allowing to encompass some discontinuous functions), \cite{chotard2019verifiable} establishes the $\varphi$-irreducibility of \eqref{eq:def-MC} based on the stability of the associated deterministic control model. It indeed extends the results of \cite{meyn1989stochastic,meyn1991asymptotic,meyn2012markov} which cover the case where  $F$ is infinitely differentiable and $\alpha(x,u)=u$. On the other hand, the results in \cite{chotard2019verifiable} can be applied to show the $\varphi$-irreducibility of Markov chains following models \eqref{eq:NLstate-space-model} relaxing smoothness conditions on $G$. Finally, compared with \cite{meyn2012markov}, the reference \cite{chotard2019verifiable} introduces the notion of steadily attracting states, which simplifies the characterization of the aperiodicity of the model.}


A particularly relevant algorithm of the form \eqref{eq:def-MC} in Evolution Strategies (ES) is ES with Covariance Matrix Adaptation (CMA-ES) \cite{hansen2001completely,hansen2003reducing} often regarded as the state-of-the-art algorithm for numerical derivative-free optimization of difficult problems 
with tremendous applications in many domains (e.g., in biology \cite{bieler2014robust,rodriguez2006hybrid}, medicine \cite{patte2022estimation}, 
machine learning \cite{akiba2019optuna,ha2018recurrent})\footnote{As of September 2023, the two main Python implementations of the CMA-ES algorithm \href{https://www.pepy.tech/projects/cma}{cma} and \href{https://www.pepy.tech/projects/cmaes}{cmaes} have more than 5 millions and 45 millions downloads respectively.}. 
Yet, while we have ample empirical evidences of its linear convergence on wide classes of functions, 
a convergence proof together with a convergence rate is still an open question. 
In order to extend linear convergence results from step-size adaptive ES \cite{auger2016linear,toure2023global} to CMA-ES, 
a first step is to show the irreducibility and topological properties of the  kernel associated to a normalized Markov chain underlying the algorithm. 
However, previous works \cite{meyn1991asymptotic,chotard2019verifiable} cannot be applied since (i) the state-space $\cX$ of this chain is a smooth manifold whereas previous analysis supposed that they were open subsets of a Euclidean space, (ii) the function $F$ is supposed to be continuously differentiable in existing results
while certain step-size updates used in CMA-ES are only locally Lipschitz.   
One motivation of the present paper is to resolve these two limitations and pave the way to a complete convergence analysis of CMA-ES.

{ In this context, the objective of this paper is to extend the theory developed in \cite{meyn1989stochastic,meyn1991asymptotic,meyn2012markov} and further expanded in \cite{chotard2019verifiable} in two directions: (1) by allowing $\cX$ and $\cW$ to be smooth manifolds rather than open subsets of Euclidean spaces, and (2) by assuming that $F$ in \eqref{eq:def-MC} is only locally Lipschitz instead of continuously differentiable.
	Under these new assumptions, analyzing the stability of \eqref{eq:def-MC} and its control requires additional arguments and new tools, which we develop here. In particular, Appendix~\ref{sec:appendix-clarke} adapts Clarke’s derivative to the setting of smooth manifolds.}

The paper is organized as follows.
In  \Cref{sec:model}, we provide a precise definition of the class of nonlinear state-space models under investigation.
In \Cref{sec:assumptions}, we outline the assumptions necessary for establishing our main results and deriving the irreducibility and aperiodicity of our model. Our main results are presented in  \Cref{sec:main-results} and are subsequently applied in  \Cref{sec:application} to
{an auto-regressive Riemannian model in \Cref{sec:random-walk} and to}
two zeroth-order optimization algorithms in \Cref{sec:cmaes,sec:csaes}.
Finally, proofs are gathered in \Cref{sec:proofs}.
Note that some of the proofs and useful definitions are given in the appendix.

\section{Main results} 

\subsection{The model and assumptions} \label{sec:model}
Let $\msx, \msw$ be two {(smooth, connected)} manifolds (see
\Cref{def:manifolds}) of dimensions $n$ and $p$ respectively,
endowed with
their Borel $\sigma$-fields denoted by $\mcbb(\msx)$ and
$\mcbb(\msw)$ respectively.  We let $\distX$ and $\distW$ be
two distance functions on $\msx$ and $\msw$ which induce the topology of 
$\msx$ and $\msw$ respectively. As a consequence of \cite[Proposition 13.2, Theorem 13.29]{lee2012introduction}, such distance functions always exist.

We consider in this paper Markov chains taking values in
$\msx$ and associated with the general recursion
\eqref{eq:def-MC}. {Throughout the paper, we denote by $P$ the
	Markov kernel associated to \eqref{eq:def-MC}.} As
emphasized in the introduction, this class of models is a
natural extension of nonlinear state-space models defined on
manifolds.

{
	As an illustration, we consider a simple example: functional Riemannian random walk models. Here, $\cX$ is assumed to be a smooth Riemannian manifold, and we denote by $(x,u) \in T \msx \mapsto \Exp_x(u)$ the exponential map on the tangent bundle $T\msx$ (see \cite[Chapter 5]{lee1997riemannian}). For clarity, we further assume that $\msx$ is complete, simply connected, and has nonpositive sectional curvature; such a manifold is called a Hadamard manifold. This class of manifolds has been extensively studied in the optimization literature; see, e.g., \cite{li2009monotone,bacak2014convex}. This assumption ensures, by the Hadamard theorem \cite[Theorem 3.1 Chapter 7]{do1992riemannian}, the existence of a global frame, i.e., a map $\iota : (x,u) \in \msx \times \rset^n \mapsto \iota_x(u) \in T\msx$ which is a smooth diffeomorphism. Therefore, without loss of generality, we identify $T\msx$ with $\msx\times \rset^n$ and regard $\Exp$ as a map from $\msx\times \rset^n$ to $\msx$. In this context, define $\{\phi_k\}_{k\in\nset}$ by the recursion
	\begin{equation}
		\label{eq:generalized-ula}
		\phi_{k+1} = \Exp_{\phi_k}(-\gamma s(\phi_k) +  U_{k+1})\eqsp,
	\end{equation}
	where $\gamma>0$ is a fixed step-size, $s: \msx \to T\msx$ is a vector field, and $\{U_{k+1}\}_{k\in\nset}$ is an \iid\ process on $\rset^n$. 
	The update equation associated with $\{\phi_{k}\}_{k\in\bN}$ can then be written in the form of \eqref{eq:def-MC} with $\msw = \rset^n $, $F(x,w)=\Exp_x(w)$ and $\alpha(x,u)=-\gamma s(x) + u$. When $s$ is the Riemannian gradient of a potential function $f:\msx\to\rset$ and $U_{k+1}/\sqrt{\gamma} \sim \mathcal{N}(0,I_n)$, \eqref{eq:generalized-ula} corresponds to the Riemannian Langevin Monte Carlo method studied in \cite{cheng2022efficient,li2023riemannian} for sampling from a distribution on $\msx$ with density proportional to $x \mapsto \exp(-f(x))$. Finally, natural extension consists in replacing the exponential map in \eqref{eq:generalized-ula} by retraction maps \cite{absil2012projection,bharath2025sampling}.}

{To further illustrate the relevance of models \eqref{eq:def-MC}, we consider in \Cref{sec:application} additional examples arising from evolution strategies (ES) for zeroth-order optimization methods~\cite{rechenberg1973evolutionsstrategie}. In particular, we apply the theory developed in this section to analyze a simplified variant of CMA-ES~\cite{hansen2001completely}, as well as a step-size adaptive ES that uses a nonsmooth step-size update.}

\subsubsection{Assumptions} \label{sec:assumptions}

We consider the following assumptions on the functions $F$ and $\alpha$ to establish ergodicity of the Markov kernel defined via \eqref{eq:def-MC}:
\begin{assumptionH}
	\label{ass:alpha}
	For any $x\in\cX$, the distribution $\mu_x$ of the random variable $\alpha(x,U_1)$ admits a density, denoted by $p_x$, with respect to a $\sigma$-finite measure $\zeta_{\cW}$, such that:
	\begin{enumerate}[label = (\roman*)]
		\item The function $(x,w)\mapsto p_x(w)$ is lower semicontinuous (l.s.c.), i.e., for any $(\bar x,\bar w)\in\cX\times\cW$, $\liminf_{(x,w)\to(\bar x,\bar w)} p_x(w)\geqslant p_{\bar x}(\bar w)$.
		\item For any $\msa \subset \mcbb(\msw)$, $\zeta_{\cW}(\msa)=0$ if and only if $\msa$ is negligeable, i.e., $\Leb(\varphi(\msa\cap U))=0$ for any chart $(\varphi,U)$ of $\msw$, where $\Leb$ stands for the Lebesgue measure.
	\end{enumerate}
\end{assumptionH}
The condition \Cref{ass:alpha} is a generalization {of \cite[A4]{meyn1991asymptotic} and} \cite[A4]{chotard2019verifiable}, where $\cW\subset\bR^p$ was instead an open subset of an Euclidean space and $\zeta_{\cW}$ the Lebesgue measure. 
If $\cW$ is equipped with a smooth Riemannian metric which makes $\cW$ a Riemannian manifold, 
a $\sigma$-finite measure satisfying \Cref{ass:alpha}(ii) would be the Lebesgue-Riemann volume measure~\cite[Chapter XII and Proposition XII.1.6]{amann2009analysis}.

We assume moreover the following on the map $F\colon\cX\times\cW\to\cX$.
\begin{assumptionH}
	\label{ass:F}
	The map $F\colon\cX\times\cW\to\cX$ is locally Lipschitz, see \Cref{def:locallyLipschitz}, on $\msx \times \msw$ with respect to the distance $\distX\oplus\distW$, defined by $\distX\oplus\distW((x,w),(x',w'))=\distX(x,x')+\distW(w,w')$ for every $((x,w),(x',w'))\in(\cX\times\cW)^2$. 
\end{assumptionH}
This assumption encompasses the requirement that $F$ be infinitely differentiable in \cite{meyn1991asymptotic,meyn2012markov}, as well as the condition of continuous differentiability considered in \cite{chotard2019verifiable}, in the case where $\cX$ and $\cW$ are open subsets of some Euclidean spaces. Indeed, any continuously differentiable or infinitely differentiable function is in particular locally Lipschitz.

For our last assumption regarding the functions $F$ and $\alpha$, we need to introduce further notations and notions {introduced in \cite{meyn2012markov,chotard2019verifiable}}.
The \textit{extended transition map} $S_x^k\colon\cW^k\to\cX$ can be defined inductively via
\begin{align}
	\label{eq:ext-transition-map}
	S_x^{k+1} (w_{1:k+1}) &\coloneqq F(S_x^k(w_{1:k}),w_{k+1}) \eqsp, \quad
	S_x^0 \coloneqq x \eqsp,
\end{align}
for $k\in\bN$, $x\in\cX$ and $w_{1:k+1}=(w_1,\dots,w_{k+1})\in\cW^{k+1}$. {The value $S_x^k(w_{1:k})$ corresponds to the $k^\mathrm{th}$ iterate of the chain $\phi_k$ defined via \eqref{eq:def-MC}, conditionally to $\phi_0=x$ and $\alpha(\phi_t,U_{t+1})=w_t$ for $t=0,\dots,k-1$.}
Remark that, by composition, if $F$ is continuous(ly locally Lipschitz), then so is $(x,w_{1:k})\mapsto S_x^k(w_{1:k})$.
Similarly, we define the extended probability density $p_x^k$ via
\begin{align}
	\label{eq:ext-proba-density}
	p_x^{k+1} (w_{1:k+1}) &\coloneqq p_x^k(w_{1:k})p_{S_x^k(w_{1:k})}(w_{k+1}) , \quad
	p_x^1 (w_1)  \coloneqq p_x(w_1) \eqsp.
\end{align}
The function $p_x^k$ is then the density of the random variable $(\alpha(\phi_0,U_1),\dots,\alpha(\phi_{k-1},U_k))$, 
with $\phi_0=x$, w.r.t.\ the product measure $\zeta_{\cW}^{\otimes k}$. 
If $(x,w)\mapsto p_x(w)$ is l.s.c., then $(x,w_{1:k})\mapsto p_x^k(w_{1:k})$ is l.s.c.\ as well. 
In this case, the control sets
\begin{equation}
	\cO_x^k \coloneqq \{ w_{1:k}\in\cW^k \mid p_x^k(w_{1:k})>0\}
\end{equation}
are nonempty open subsets of $\cW^k$. 
The control set $\cO_x^k$ corresponds to the set of paths $w_{1:k}$ starting at $x$ which have positive density $p_x^k(w_{1:k})$.

\begin{figure}
	\hspace*{2cm}\includegraphics[scale=0.2]{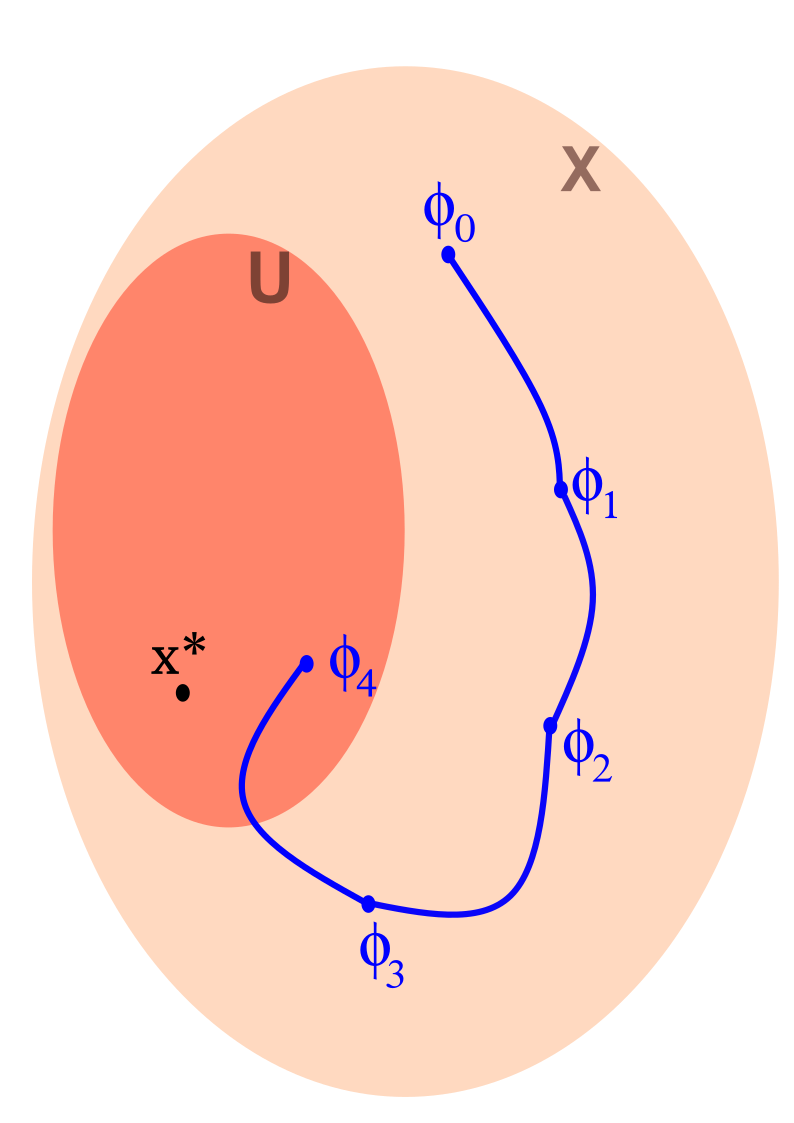} \hfill
	\includegraphics[scale=0.2]{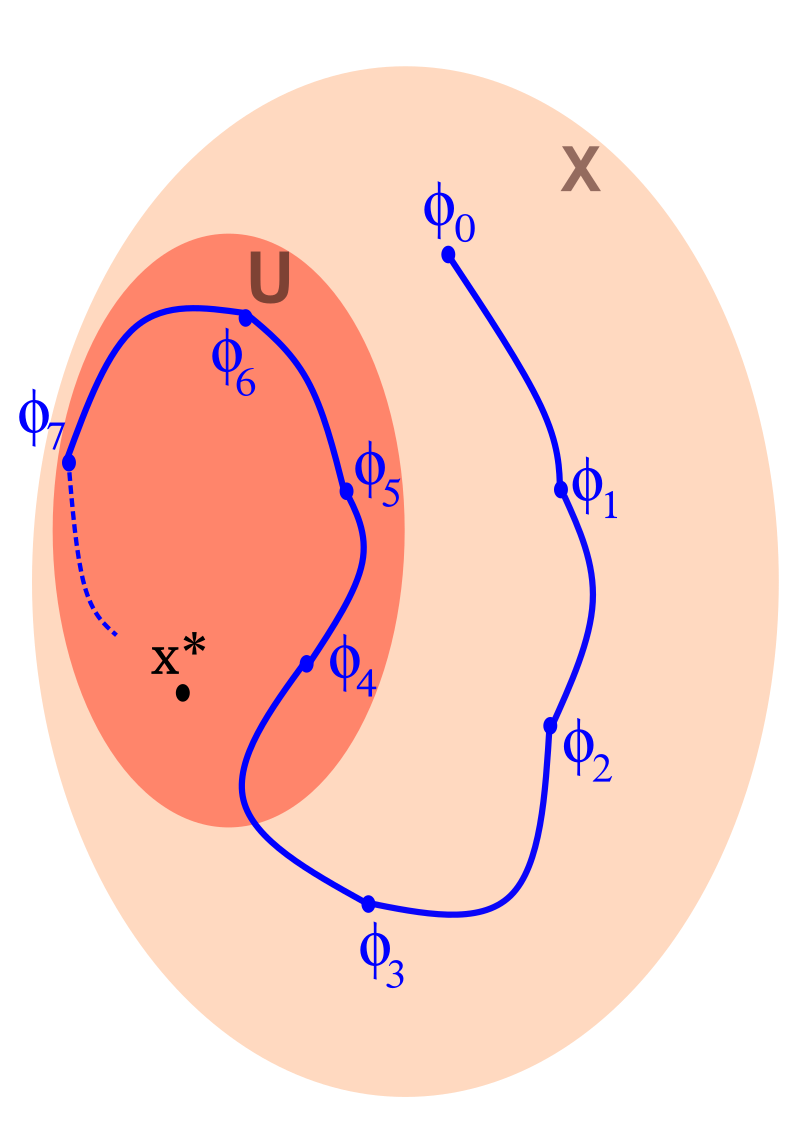}\hspace*{2cm}
	\caption{Left: Illustration of a globally attracting state $x^*$, for any neighborhood $U$ of $x^*$ and any starting state $\phi_0$, there exists a $k$-steps path from $\phi_0$ to $U$. \\		
		Right: Illustration of a steadily attracting state $x^*$, for any neighborhood $U$ of $x^*$ and any starting state $\phi_0$, there exist $T>0$ and $k$-steps paths from $\phi_0$ to $U$ for every $k\geqslant T$.}
	\label{fig:attracting-states}
\end{figure}

Moreover, for $x\in\cX$, $A$ a measurable subset of $\cX$, and $k>0$, we say that $w_{1:k}\in\cW^k$ is a $k$\textit{-steps path} from $x$ to $A$ if $w_{1:k}\in\cO_x^k$ and $S_x^k(w_{1:k})\in A$, 
implying that $A$ is then reachable by $P$ from $x$. 
A point $x^*\in\cX$ is said to be a \textit{globally attracting state} if for any $y\in\cX$ and any neighborhood $U$ of $x^*$, there exist $k>0$ and a $k$-steps path between $y$ and $U$ (the original definition of a globally attracting state is actually given in \eqref{eq:def-globally-attracting-state} and we show in Proposition~\ref{p:characterization-globally-attracting} the equivalence with this latter statement). It is said to be \textit{steadily attracting} if for any $y\in\cX$ and any neighborhood $U$ of $x^*$, there exists $T>0$ such that for every $k\geqslant T$, we can find a $k$-steps path between $y$ and $U$. Note that any steadily attracting state is in particular globally attracting. 
These two notions are illustrated in \Cref{fig:attracting-states}.
{The notion of globally attracting states was introduced in \cite{meyn1991asymptotic,meyn2012markov}. Their characterization through $k$-step paths was later given in \cite{chotard2019verifiable}, which also introduced the concept of steadily attracting states.}

{Assuming \Cref{ass:alpha} and \Cref{ass:F}, as emphasized in \Cref{t:sufficient-condition-ireducibility}, we show that the kernel $P$ defined via \eqref{eq:def-MC} is $\varphi$-irreducible exhibiting the existence of a globally attracting state (it is in fact an equivalence)
	. On a related note, we deduce in \Cref{t:sufficient-condition-aperiodicity} that the existence of a steadily attracting state is equivalent to the $\varphi$-irreducibility and aperiodicity of $P$.

\newcommand{\controllability}[1]{\textnormal{(}\hyperref[eq:controllabilty-condition]{\ensuremath{\mathrm{C}_{#1}}}\textnormal{)}}

We introduce now the notation $\partial f$ for the \textit{Clarke's generalized Jacobian} of a locally Lipschitz function $f\colon\cX\to\cY$ between two\del{ metrized} manifolds $\cX$ and $\cY$. These Jacobians have been defined in \cite{clarke1990optimization}, and we recall the definition in the Euclidean case in \Cref{def:clarke-jacobian}. For the sake of completeness, we define here and give basic properties in \Cref{sec:appendix-clarke} of the Clarke's Jacobian for functions defined on manifolds.

\begin{definitionProposition}[Clarke's generalized Jacobian on manifolds]
	\label{p:clarke-derivative-manifold-chain-rule}
	Let $\cX$ and $\cY$ be two manifolds and  $f\colon \cX\to\cY$ be locally Lipschitz at $x_0\in\cX$. Let $(\varphi,U)$ be a local chart of $\cX$ around $x_0$ and $(\psi,V)$ be a local chart of $\cY$ around $f(x_0)$. Define $g=\psi\circ f\circ \varphi^{-1}$. Then $g\colon\varphi(U)\to\psi(V)$ is locally Lipschitz at $\varphi(x_0)$, and we can define
	\begin{equation}
		\label{eq:def-jac-manifolds}
		\partial f (x_0) = \left\{ \cD\psi^{-1}(g\circ\varphi(x_0))\circ h \circ \cD \varphi (x_0) \mid h\in\partial g (\varphi(x_0)) \right\} ,
	\end{equation}
	where $\cD$ denotes the usual differential operator, and $\partial$ the Clarke differential operator.
	This definition does not depend on the choice of the charts $(\varphi,U)$ and $(\psi,V)$.
\end{definitionProposition}
\begin{proof}
	See \Cref{sec:appendix-clarke}.		
\end{proof}
In the case of a differentiable function $f$, the definition of Clarke's generalized Jacobian corresponds to the definition of the Jacobian, i.e., $\partial f(x)=\{\cD f(x)\}$.
{The notion of Clarke's generalized Jacobian is used to formulate the controllability condition for an element $x\in\cX$: 
	\begin{equation}
		\tag{\ensuremath{\mathrm{C}_x}}
		\label{eq:controllabilty-condition}
		\text{\del{T}there exists } w_{1:k}\in\overline{\cO_{x}^k} \text{ such that } \partial_w S_x^k({w}_{1:k}) \text{ is of maximal rank.}
	\end{equation} 
	Note that here, $\partial_w S_x^k({w}_{1:k})$ is of maximal rank, is understood as any element of {the Clarke's derivative with respect to $w_{1:k}$} $\partial_w S_x^k(w_{1:k})$ is of rank $n$, the dimension of $\cX$.
	In \Cref{c:Tchain}, we show that $P$ is a T-chain assuming that condition \controllability{x} holds for every state $x\in\cX$.
	In comparison to \Cref{ass:sub_differe_S} below, we do not assume that states $x$ for which \controllability{x} holds are globally attracting. However, we show that if \controllability{x^*} holds for $x^*$ a globally attracting state, then it holds for every state in $\cX$, see \Cref{p:steadily-att-state-differentiable}.
	
	\begin{assumptionH}
		\label{ass:sub_differe_S}
		The controllability condition \controllability{x^*} is satisfied for a globally attracting state $x^*$.
	\end{assumptionH}
	
	Alternatively, if we want to prove aperidocity on top of $\varphi$-irreducibility, we assume instead the following.
	\begin{assumptionH}
		\label{ass:sub_differe_S_steadily}
		The controllability condition \controllability{x^*} is satisfied for a steadily attracting state $x^*$.
	\end{assumptionH}
	{Remark that \Cref{ass:sub_differe_S_steadily} implies \Cref{ass:sub_differe_S}.} 
	Assumptions \Cref{ass:sub_differe_S}-\Cref{ass:sub_differe_S_steadily} also appear in \cite{chotard2019verifiable} but with the additional condition that the functions $S_x^k$ are continuously differentiable for $x\in\cX$ and $k>0$, condition that we relax here.
	{Condition \Cref{ass:sub_differe_S} was first introduced in \cite{meyn1991asymptotic,meyn2012markov}, while \Cref{ass:sub_differe_S_steadily} was later considered in \cite{chotard2019verifiable}.}
	{Globally and steadily attracting states are characterized by \Cref{p:characterization-globally-attracting} and \Cref{p:sufficient-condition-steadily-attracting}(ii) respectively below.
		In \Cref{sec:application}, we give one example of a smooth model on manifolds, and one example of a nonsmooth model on a Euclidean space, for which we show that \Cref{ass:sub_differe_S_steadily} holds.}

	\subsection{Main results}\label{sec:main-results}
	Before stating our main results, we introduce concepts that are needed for their statements.
	Given $P$ a Markov kernel on $(\cX,\cB(\cX))$, we define $P^1=P$ and for $k\geqslant1$, $x\in\cX$ and $\mathsf{A}\in\cB(\cX)$, $P^{k+1}(x,\mathsf{A})= \int P(y,\mathsf{A})P^k(x,\mathrm{d}y)$.
	We say that $P$ is $\varphi$-irreducible when there exists a nontrivial measure $\varphi$ on $\cB(\cX)$ such that for any $\msa\in\cB(\cX)$ with $\varphi(\msa)>0$, we have
	$
	\sum_{k\geqslant1}P^k(x,\msa)>0 \text{ for every } x\in \cX .
	$
	Let $b$ be a probability distribution on $\bN$, and let $K_b$ be the transition kernel defined by 
	$
	K_b(x,\msa) \coloneqq \sum_{k\geqslant0} b(k) P^k(x,\msa) .
	$
	A substochastic transition kernel $T$ with $K_b\geqslant T$
	such that $x\mapsto T(x,\msa)$ is lower semicontinuous for every $\msa\in\cB(\cX)$ is called a \textit{continuous component} of $K_b$. If $P$ admits a distribution $b$ such that there exists a continuous component $T$ of $K_b$ with $T(\cdot,\cX)>0$, then $P$ is called a \textit{T-chain}.
	
	A set $C\in\cB(\cX)$ is called \textit{petite} if there exist a probability distribution $b$ on $\bN$ and a nontrivial measure $\nu_b$ on $\cB(\cX)$ such that
	$
	K_b(x,\msa)\geqslant\nu_b(\msa) $ for every $x\in\cX$ and $\msa\in\cB(\cX)$. If moreover $b=\delta_a$ the Dirac distribution at some $a\in\bN$, then $C$ is called $a$-small.

	If $P$ is $\varphi$-irreducible, then the family $(D_i)_{i=1,\dots,d}\in\cB(\cX)^d$ is called a \textit{$d$-cycle} when
	\begin{equation}
		\label{eq:aperiodicity-def}
		\left\{
		\begin{array}{l}
			P(x,D_{i+1})=1 \text{ for } x\in D_i \text{ and } i=0,\dots,d-1 ~\mathrm{mod}~d \\
			\varphi ((\cup_{1\leqslant i\leqslant d} D_i)^c) = 0 \mbox{ for any irreducibility measure } \varphi \eqsp.
		\end{array}
		\right.
	\end{equation}
	By \cite[Theorem 5.4.4 and Proposition 5.2.4]{meyn2012markov}, if $P$ is $\varphi$-irreducible, then there exist $d\geqslant1$ and a $d$-cycle. The \textit{period} of $P$ is the largest integer $d$ for which there exists a $d$-cycle. If the period of $P$ is equal to $1$, then $P$ is said to be \textit{aperiodic}.

	We have now all the tools to state our main contribution. 
	\begin{theorem}
		\label{theo:1}
		Assume \Cref{ass:alpha}-\Cref{ass:F} and \Cref{ass:sub_differe_S}. Then, the
		Markov kernel $P$ defined via \eqref{eq:def-MC} is a $\varphi$-irreducible T-chain, and any compact set is petite. 
		If moreover \Cref{ass:sub_differe_S_steadily} holds, then, the
		Markov kernel $P$ is aperiodic, and any compact set is small. 
	\end{theorem}
}%
In
addition to the assumptions of \Cref{theo:1}, if $P$ is positive recurrent {(i.e., $P$ is $\varphi$-irreducible and admits an invariant probability measure)}, then $\seq{\phi_{\t}}{\t \in \nset}$ is ergodic, i.e., $P$ admits a unique stationary distribution $\pi$ and for $\pi$-almost every $x\in\msx$,
\begin{equation}
	\label{eq:1}
	\lim_{k \to \plusinfty} \tvnorm{\updelta_x P^k - \pi} = 0 \eqsp. 
\end{equation}
Moreover, if we suppose that $P$ is Harris recurrent, then a Law of Large Numbers holds, see~\cite[Theorem 17.0.1]{meyn2012markov}. For any $\pi$-integrable function $g$, a Markov chain $\{\phi_k\}_{k\in\bN}$ associated to the kernel $P$ satisifies
\begin{equation}
	\lim_{T\to\infty} \frac{1}{T} \sum_{k=0}^{T-1} g\left(\phi_k\right) = \int g\, \mathrm{d}\pi \eqsp. 
\end{equation}
{Harris recurrence can be established typically as a consequence of a Foster-Lyapunov condition~\cite[Theorem 13.0.1]{meyn2012markov}, i.e., function $V\colon\cX\to [0,+\infty]$ finite at least at one point of $\cX$, of a petite set $C$ and of a constant $b<\infty$, such that, for any $x\in\cX$, we have
	\begin{equation}
		\label{eq:foster-lyapunov}
		\int V(y) P(x,\mathrm{d}y) - V(x) \leqslant -1 + b \1\{ x\in C\} \eqsp.
	\end{equation}
}%
The proof of \Cref{theo:1} is postponed to \Cref{sec:proofs-main-results}. 
It relies on intermediary results that to a great extent are generalizations of results in
{\cite{meyn1991asymptotic,meyn2012markov,chotard2019verifiable}} {when we assume \Cref{ass:sub_differe_S_steadily}}.
In particular, \Cref{p:characterization-globally-attracting,p:equivalence-ksteps,p:sufficient-condition-steadily-attracting} characterize globally attracting states, 
reachable states and steadily attracting states respectively. 
\Cref{p:rank-condition-on-all-states,p:steadily-att-state-differentiable,p:implicit-Skx,p:forward-accessibility} provide consequences of the assumption of controllability \eqref{eq:controllabilty-condition}. 
\Cref{l:small-set-lemma} is a generalization of \cite[Lemma 3.0]{meyn1991asymptotic}, 
which turns out to be  useful to prove that the controllability condition \eqref{eq:controllabilty-condition} implies that the Markov kernel $P$ is a T-chain, 
as stated in \Cref{p:smallset-property} and \Cref{c:Tchain}. 
\Cref{p:support-irreducibility-measure} characterizes the support of the irreducibility measures of $P$,
while \Cref{t:sufficient-condition-ireducibility,t:practical-condition-ireducibility,t:sufficient-condition-aperiodicity,t:practical-condition-aperiodicity} end 
the proof of \Cref{theo:1}.

\section{Applications}~
\label{sec:application}

{\input{auto-reg}}

\subsection{An instructive example: CMA-ES}\label{sec:cmaes}
We introduce here a simplified version of the {numerical optimization} algorithm {called} evolution strategy with covariance matrix adaptation
(CMA-ES)\del{,  first developed in}
\cite{hansen2001completely,hansen2003reducing}, \del{with}{which, for an objective function $f: \bR^d \to \bR$, aims}\del{ the aim} to 
solve:
\begin{equation}
	\tag{P}\label{eq:(P)}
	\text{find } x^*\in\underset{x\in\bR^d}{\Argmin}~ f(x) \eqsp.
\end{equation}
To this end, it approximates the optimum $x^*$ of the {objective} function $f$ by a multivariate normal distribution $\mathcal{N}(m_{\t},C_{\t})$ for a mean $m_{\t}\in\bR^d$ and a covariance matrix $C_{\t}\in\Sdpp$ that are updated iteratively. More precisely, for each $\t\in\bN$, given $m_{\t}\in\bR^d$ and $C_{\t}\in\Sdpp$, the algorithm can be described as follows. First, a population of $\lambda\geqslant2$ offspring is sampled using
\begin{equation}\label{eq:sampling-cma}
	U_{\t+1}^1,\dots,U_{\t+1}^\lambda \sim \mathcal{N}(0,I_d) \quad \text{i.i.d.\ and independently of } (m_\t,C_\t)\eqsp,
\end{equation}
so that, conditionally to $(m_k,C_k)$, the offspring satisfy $m_k+\sqrt{C_k}U_{k+1}^i\sim\mathcal{N}(m_k,C_k)$, for $i=1,\dots,\lambda$.
Next, we rank the offspring so that we define a permutation $s_{k+1}\in\mathfrak{S}_{\lambda}$ of $\{1,\dots,\lambda\}$ satisfying
\begin{equation}\label{eq:ranking-cma}
	f\left(m_{\t}+\sqrt{C_{\t}}U_{\t+1}^{s_{\t+1}(1)}\right)\leqslant \dots \leqslant f\left(m_{\t}+\sqrt{C_{\t}}U_{\t+1}^{s_{\t+1}(\lambda)}\right) \eqsp.
\end{equation}
Then{, given the $\mu\in\{1,\dots,\lambda\}$ best offspring}, \del{we update }the mean {is moved towards the best solutions with the following update}\del{ according to}
\begin{equation}\label{eq:update-mean-cma}
	m_{\t+1} = m_{\t} + \sqrt{C_{\t}} \sum_{i=1}^\mu w_i U_{\t+1}^{s_{\t+1}(i)} \eqsp,
\end{equation}
and the covariance matrix update reads
\begin{equation}
	C_{\t+1} = (1-c) C_{\t} + c \sqrt{C_{\t}}\left(\sum_{i=1}^\mu w_i \left(U_{\t+1}^{s_{\t+1}(i)}\right)\left(U_{\t+1}^{s_{\t+1}(i)}\right)^\top \right)\sqrt{C_{\t}} \eqsp.
\end{equation}
It increases the likelihood to sample in the directions where good solutions were found.
In the above equations, the weights $w_1\geqslant \dots\geqslant w_\mu >0$ satisfy $\sum_{i=1}^\mu w_i=1$, and we call $c\in(0,1)$ the learning rate  for the covariance matrix.
In ES, the function values are not used explicitly to update the state variables. 
It influences the update only through the ranking of candidate solutions via the permutation $s_{k+1}$. Consequently, the algorithms are invariant with respect to strictly increasing transformations of the objective function (that preserve the ranking). 
In this context, a natural class of functions to analyze the convergence of ES are scaling-invariant functions \cite{auger2016linear,toure2021scaling}. A function $f$ is said to be scaling-invariant w.r.t.\ $x^*$ if, for every $x,y\in\bR^d$ and $\rho>0$, we have
\begin{equation}
	\label{eq:scaling-invariant}
	f(x+x^*) \leqslant f(y+x^*) \Leftrightarrow f(\rho x+x^*)\leqslant f(\rho y+x^*) .
\end{equation}
{Convergence of step-size adaptive ES on scaling-invariant functions with smooth level sets was established --for specific assumptions on the step-size update-- in previous work \cite{toure2023global}.}
\del{Scaling-invariant functions have been of primary interest in optimization \del{through}{with} evolution strategies, see \cite{toure2021scaling}. }
Assuming that the objective function $f$ satisfies \eqref{eq:scaling-invariant}, we define then the following quantities
\begin{equation}\label{eq:normMC}
	z_{\t} = \frac{m_{\t}-x^*}{\sqrt{R(C_{\t})}} \quad ; \quad \ncovmat_{\t} = \frac{C_{\t}}{R(C_{\t})}\eqsp, 
\end{equation}
where $R=\det(\cdot)^{1/d}\colon\Sdpp\to\bR_+$. {We assume w.l.o.g.\ that $x^*=0$.}
Then the sequence $\seq{(z_{\t},\ncovmat_{\t})}{\t\in\bN}$ defines a time-homogeneous Markov chain
which obeys to the model \eqref{eq:def-MC}, see \Cref{p:cma-control-model}, with $\cX=\bR^d\times R^{-1}(\{1\})$, $\cU=\bR^{d\times\lambda}$, $\cW=\bR^{d\times\mu}$, and
\begin{equation}
	\label{eq:FCMA}
	\begin{array}{rl}
		F\colon & \cX\times\cW \to\cX \\
		& \left((z,\ncovmat),(v_1,\dots,v_\mu)\right)  \mapsto \left( \cfrac{z+\sqrt{\ncovmat}\sum_{i=1}^\mu w_i v_i}{R^{1/2}\left(K(\ncovmat,v_1,\dots,v_\mu)\right)} , \cfrac{K(\ncovmat,v_1,\dots,v_\mu)}{R\left(K(\ncovmat,v_1,\dots,v_\mu)\right)} \right)
		\eqsp,
	\end{array}
\end{equation}
where
$$
K(\ncovmat,v_1,\dots,v_\mu) = (1-c) \ncovmat+c \sqrt{\ncovmat} \left( \sum_{i=1}^\mu w_i v_iv_i^\top \right)\sqrt{\ncovmat} \eqsp,
$$
and with
\begin{equation}
	\label{eq:alphaCMA}
	\begin{array}{rl}
		\alpha\colon & \cX\times\cU\to\cW \\
		& \left((z,\ncovmat),(u_1,\dots,u_\lambda)\right)  \mapsto \left(u_{s(1;z,\ncovmat,u_{1:\lambda})},\dots,u_{s(\mu;z,\ncovmat,u_{1:\lambda})}\right) \eqsp,
	\end{array}
\end{equation}
where given $u_{1:\lambda}=(u_1,\ldots,u_\lambda)\in(\bR^d)^\lambda$, $z\in\bR^d$ and $\ncovmat\in\Sdpp$, we denote by $s(\cdot;z,\ncovmat,u_{1:\lambda})$ a permutation that sorts the $f(z+\sqrt{\ncovmat}u_i)$, $i=1,\dots,\lambda$. To ensure uniqueness of this permutation, we impose a tie-break, e.g., if $i<j$ are such that $f(z+\sqrt{\ncovmat}u_i)=f(z+\sqrt{\ncovmat}u_j)$, then $s(\cdot;z,\ncovmat,u_{1:\lambda})^{-1}(i)<s(\cdot;z,\ncovmat,u_{1:\lambda})^{-1}(j)$. {Note that $\cX$ is not an open subset of a Euclidean space, hence the results in \cite{chotard2019verifiable} do not apply and neither results in \cite{meyn1991asymptotic,meyn2012markov}. However,  $\msx$ is a smooth manifold by the preimage theorem, see e.g., \cite[Chapter 1, Section 4]{guillemin2010differential}). 
}
We show in \Cref{sec:application} that our results apply and we prove that $\seq{(z_{\t},\ncovmat_{\t})}{\t\in\bN}$ defines a $\varphi$-irreducible aperiodic T-chain and that all compact subsets of $\cX$ are small.


If we establish moreover that the chain $\seq{(z_\t,\ncovmat_\t)}{\t\in\bN}$ is positive recurrent, then we obtain that CMA-ES behaves linearly, as stated below.

\begin{theorem}
	{Consider a scaling-invariant function with respect to $x^*$ and the Markov chain $\seq{(z_\t,\ncovmat_\t)}{\t\in\bN}$ defined in \eqref{eq:normMC} ensuing  from CMA-ES minimizing $f$.}
	Suppose that \del{the Markov chain }$\seq{(z_\t,\ncovmat_\t)}{\t\in\bN}$ is {a} $\varphi$-irreducible aperiodic positive recurrent chain with invariant probability measure $\pi$. If the function $(z,\ncovmat)\mapsto\log\|z\|$ is $\pi$-integrable on $\bR^d\times\Sdpp$, then almost surely we have
	\begin{equation}
		\label{eq:linear-convergence}
		\lim_{\t\to\infty} \frac{1}{\t}\log\frac{\|m_\t-x^*\|}{\|m_0-x^*\|} = -\mathrm{CR} \in\bR \eqsp.
	\end{equation}
	When moreover $\mathrm{CR}>0$, we say that CMA-ES converges linearly to $x^*$.
\end{theorem}

\begin{proof}
	{
		Assume that $x^*=0$. Since CMA-ES is invariant by translation~\cite{auger2016these}, \eqref{eq:linear-convergence} would generalize to any value of $x^*$.
	}
	Since $\seq{(z_\t,\ncovmat_\t)}{\t\in\bN}$ is supposed to be $\varphi$-irreducible, aperiodic and positive recurrent, by \cite[Theorem 17.0.1]{meyn2012markov}, we know that for all $\pi$-integrable function $g$, we have that 
	\begin{equation}\label{eq:LLN}
		\lim_{T\to\infty} \frac{1}{T}\sum_{\t=0}^{T-1} g(z_\t,\ncovmat_\t) = \int g(z,\ncovmat) \d\pi(z,\ncovmat).
	\end{equation} 
	However, we have
	\begin{align}
		\frac{1}{T} & \nonumber \log \frac{\| m_T \|}{\| m_0  \|} = \frac{1}{T} \sum_{\t=0}^{T-1} \left( \log\| m_{\t+1}  \| -\log\| m_{\t} \| \right) \\\label{eq:FLLN1}
		=& \frac{1}{T} \sum_{\t=0}^{T-1} \left( \log\| z_{\t+1} \| -\log\| z_\t\| \right) +  \\& \frac{1}{2dT}  \sum_{\t=0}^{T-1}  \log\det 
		\underbrace{\left(
			(1-c)\ncovmat_\t+\sqrt{\ncovmat_\t}\left(\sum_{i=1}^\mu w_i \left(U_{k+1}^{s_{k+1}(i)}\right) \left(U_{k+1}^{s_{k+1}(i)}\right)^\top \right)\sqrt{\ncovmat_\t} 
			\right)}_{\eqqcolon \tilde\ncovmat_{k+1}}
		\label{eq:FLLN2} .
	\end{align}
	But, by assumption, $(z,\ncovmat)\mapsto \log\|z\|$ is $\pi$-integrable. Moreover, $\det(\ncovmat_\t)=1$, hence
	$$
	\det \left(\tilde{\ncovmat}_{\t+1}\right) = \det\left((1-c)I_d+c\sum_{i=1}^\mu w_i\left(U_{k+1}^{s_{k+1}(i)}\right)\left(U_{k+1}^{s_{k+1}(i)}\right)^\top\right) .
	$$
	Moreover,
	$$
	1-c\leqslant \det\left((1-c)I_d+c\sum_{i=1}^\mu w_i\left(U_{k+1}^{s_{k+1}(i)}\right)\left(U_{k+1}^{s_{k+1}(i)}\right)^\top\right)^{1/d} \leqslant 1-c + c \max_{i=1,\dots,\mu} \|U_{k+1}^i\|^2
	$$
	which defines an integrable quantity, since the vectors $U_{\t+1}^i$, $\t\in\bN$, $i=1,\dots,\lambda$, are standard Gaussian vectors of $\bR^d$. Applying \eqref{eq:LLN} to \eqref{eq:FLLN1} and \eqref{eq:FLLN2}, we find the stated result with 
	\begin{equation}
		\mathrm{CR} = -\frac{1}{2d} \bE_{(z,\ncovmat)\sim\pi}  \left[ \det\left( (1-c) I_d + c\sum_{i=1}^\mu w_i\left(U_{1}^{s(i;z,\ncovmat,U_1^{1:\lambda})}\right)\left(U_{1}^{s(i;z,\ncovmat,U_1^{1:\lambda})}\right)^\top  \right)    \right]  .
	\end{equation}
\end{proof}
The previous theorem illustrates how the $\varphi$-irreducibility and aperiodicity of  $\seq{(z_\t,\ncovmat_\t)}{\t\in\bN}$ are instrumental to obtain linear convergence of CMA-ES.%
\footnote{The variant of CMA-ES presented here differs significantly from the default CMA-ES (used in applications) where both step-size adaptation and covariance matrix adaptation are used. In addition, the covariance matrix update presents an additional mechanism (rank-one update). The combination of all the mechanisms is important to obtain fast convergence in many situations. This variant with however a learning rate on the mean update has been analyzed in previous theoretical works~\cite{akimoto2012theoretical}, and it has been proven to be a discretized version of a natural gradient update on the manifold of probability distributions~\cite{akimoto2010bidirectional}.}

{Let us assume} that $f$ has Lebesgue-negligible level sets, i.e., $\Leb(\cL_t)=0$, with
\begin{equation}
\cL_t \coloneqq \{x\in\bR^d \mid f(x)=t\} \quad \text{for } t\in\bR \eqsp.
\end{equation}
Stability of Markov chains defined in the context of ES with step-size adaptation has been proven~\cite{auger2016linear,toure2023global}, yielding to linear convergence. We complement these results applying now \Cref{theo:1} to show the stability of $\seq{(z_k,\ncovmat_k)}{k\in\bN}$. First, observe that the assumption \Cref{ass:F} is automatically satisfied, since $F$ is continuously differentiable. As for \Cref{ass:alpha}, we use the following result.

\begin{proposition}
\label{p:density-cma}
Suppose that $f$ has Lebesgue-negligible level sets. Define for any $\theta=(z,\ncovmat)\in\cX$ and $v=(v_1,\dots,v_\mu)\in\cW$,
\begin{multline}
	p_{\theta}(v)=\frac{\lambda!}{(\lambda-\mu)!} \1 {\left\{f\left(z+\sqrt{\ncovmat}v_1\right)<\dots<f\left(z+\sqrt{\ncovmat}v_\mu\right)\right\}} \\
	\times \left(1-Q_{\theta}^f(v_\mu)\right)^{\lambda-\mu} \gamma^d(v_1)\dots \gamma^d(v_\mu)
\end{multline}
with $Q_{\theta}^f(u)=\int\1 \{f(z+\sqrt{\ncovmat}\xi)<f(z+\sqrt{\ncovmat}u)\} \gamma^d(\xi)\d\xi$ and where $\gamma^d$ is the density of the $d$-dimensional standard normal distribution w.r.t.\ Lebesgue. Then, $p_{\theta}$ defines a density (w.r.t.\ Lebesgue in $\bR^{d\mu}$) of the random variable $\alpha(\theta,U_1)$.
\end{proposition}
If $f$ has Lebesgue-negligible level sets and is continuous, it follows that \Cref{ass:alpha} holds.

The proof of \Cref{p:density-cma} mimics the one of \cite[Proposition 5.2]{chotard2019verifiable}, but is given for completeness in \Cref{sec:additional-proofs}. Then, it remains to prove \Cref{ass:sub_differe_S_steadily} and in particular to find a steadily attracting state $\theta^*=(z^*,\ncovmat^*)$ for which there exist $k>0$ and ${v}^*_{1:k}\in\overline{\cO_{\theta^*}^k}$ such that $\partial S_{\theta^*}^k({v}^*_{1:k})$ is of maximal rank. This is achieved in the following proposition proven in \Cref{proof-for-CMA}.
\begin{proposition}
\label{p:steadily-attracting-cma}
Suppose that $f$ is continuous, scaling-invariant with Lebesgue-negligible level sets. Then, 
\begin{enumerate}
	\item[(i)] the state $\theta^*=(0,I_d)$ is steadily attracting~;
	\item[(ii)] there exists $k>0$ and ${v}^*_{1:k}\in\overline{\cO^{k}_{\theta^*}}$ such that $\cD S_{\theta^*}^{k}({v}^*_{1:k})\colon \cW^k \to T_{S_{\theta^*}^{k}({v}^*_{1:k})}\cX$ is surjective, hence is full rank, where
	$$
	T_{S_{\theta^*}^{k}({v}^*_{1:k})}\cX = \bR^d \times \ker \left(\mathcal{D}\det(I_d)\right) ,
	$$
	where $\det$ is the determinant map on the set of symmetric matrices $\Sd$, and $\ker$ denotes the kernel of a linear application.
\end{enumerate}
\end{proposition}
Then, by applying \Cref{theo:1}, the $\varphi$-irreducibility and aperiodicity of the chain $\seq{\theta_k}{k\in\bN}$ follow.
\begin{theorem}
Suppose that $f$ is continuous, scaling-invariant with Lebesgue-negligible level sets. Then the Markov chain $\seq{\theta_k}{k\in\bN}$ defines a time-homogeneous $\varphi$-irreducible aperiodic T-chain, for which any compact subset of $\cX$ is small.
\end{theorem}

\subsection{A nonsmooth example: \del{SA-ES}{a step-size adaptive ES}}\label{sec:csaes}

{We present here an other sim\-pli\-fi\-ca\-tion of CMA-ES where instead of adapting a full covariance matrix, a scaling factor called step-size is adapted such that the covariance matrix reads  $\sigma_k^2 I_d$.}
\del{	Similarly to CMA-ES, evolution strategy with step-size adaptation aims to solve the optimization problem \eqref{eq:(P)}. We approximate here t}{In this step-size adaptive algorithm, t}he optimum $x^*\in\bR^d$ {of the problem \eqref{eq:(P)} is approximated} by a multivariate normal distribution $\mathcal{N}(m_k,\sigma_k^2 I_d)$, where the mean $m_k\in\bR^d$ and the step-size $\sigma_k>0$ are updated as follows. For $k\in\bN$, given a mean $m_k\in\bR^d$ and a step-size $\sigma_k>0$, we
sample $U_{\t+1}^1,\dots,U_{\t+1}^\lambda$, rank them by defining the permutation $s_{\t+1}\in\mathfrak{S}_\lambda$ and update the mean $m_{\t+1}$ according to \eqref{eq:sampling-cma}, \eqref{eq:ranking-cma}, \eqref{eq:update-mean-cma}, respectively, where we replace $C_k$ by $\sigma_k^2 I_d$.
\del{ sample a population of $\lambda\geqslant2$ offspring {using independant samples}
\begin{equation}
	U_{k+1}^1,\dots,U_{k+1}^\lambda \sim \mathcal{N}(0,I_d) \quad \text{i.i.d.\ independently of $(m_k,\sigma_k)$,}
\end{equation}
{so that, conditionally to $(m_k,\sigma_k)$, the offspring satisfy $m_k+\sigma_k U_{k+1}^i\sim\mathcal{N}(m_k,\sigma_k^2 I_d)$, for $i=1,\dots,\lambda$.}
Then, we rank the offspring, that is, we find a permutation $s_{\t+1}$ of $\{1,\dots,\lambda\}$, so that almost surely we obtain
\begin{equation} 
	f\left( m_k+\sigma_k U_{k+1}^{s_{k+1}(1)} \right) \leqslant \dots \leqslant f\left( m_k+\sigma_k U_{k+1}^{s_{k+1}(\lambda)} \right) .
\end{equation}
Now, given the $\mu\in\{1,\dots,\lambda\}$ best offspring, {the mean is updated so as to move towards the best solutions. More precisely the new mean  equals the weighted average of the $\mu$ best solutions that writes}
\begin{equation}
	m_{k+1} = m_k + \sigma_k \sum_{i=1}^\mu w_i U_{k+1}^{s_{k+1}(i)}
\end{equation}}
The step-size update obeys
\begin{equation}
\label{eq:update-sigma}
\sigma_{k+1} = \sigma_k \times \exp\left( \frac{1}{d_\sigma} \left( \frac{\sqrt{\mueff}\left\|\sum_{i=1}^\mu w_i U_{k+1}^{s_{k+1}(i)}\right\|}{\bE\|\mathcal{N}(0,I_d)\|} -1 \right) \right)
\end{equation}
\del{We choose here the weights $w_1\geqslant\dots\geqslant w_\mu >0$ and  $\sum_{i=1}^\mu w_i=1$ as in \Cref{sec:cmaes}.}where we define $\mueff=\sum_{i=1}^{\mu}w_i^2$ and fix $d_\sigma>0$ (usually $d_\sigma\approx1$). 
Moreover, as in \Cref{sec:cmaes}, we assume $f$ to be scaling-invariant, see \eqref{eq:scaling-invariant}. W.l.o.g.\, we suppose that $f$ is scaling-invariant w.r.t.\ $x^*=0$. Then, by defining
\begin{equation}
z_k = \frac{m_k-x^*}{\sigma_k} ,
\end{equation}
we get that the sequence $\seq{z_{\t}}{\t\in\bN}$ is a time-homogeneous Markov chain which obeys to the model \eqref{eq:def-MC} {(see \cite[Proposition~4]{toure2023global})} with $\cX=\bR^d$, $\cU=\bR^{d\times\lambda}$, $\cW=\bR^{d\times\mu}$,
\begin{equation}
\label{eq:CSA-update-MC}
\begin{array}{rl}
	F\colon & \cX\times \cW \to \cX \\
	& (z,(v_1,\dots,v_\mu)) \mapsto {\left(z+\sum_{i=1}^\mu w_i v_i\right)} \times {\exp\left(-\cfrac{1}{d_\sigma}\left(\cfrac{\sqrt{\mueff}\|\sum w_iv_i\|}{\bE\|\mathcal{N}(0,I_d)\|}-1\right)\right)}
\end{array}
\end{equation}
and 
\begin{equation}
\begin{array}{rl}
	\alpha\colon & \cX\times \cU\to \cW \\
	& (z,(u_1,\dots,u_\lambda)) \mapsto  \left( u_{s(1;z,I_d,u_{1:\lambda})},\dots,u_{s(\lambda;z,I_d,u_{1:\lambda})} \right)
\end{array}
\end{equation}
where we define the permutation $s(\cdot;z,I_d,u_{1:\lambda})$ as in \Cref{sec:cmaes}.
Here, $F$ is not continuously differentiable, and we cannot use the results of \cite{chotard2019verifiable} to analyze this chain. In addition, results in \cite{meyn1991asymptotic,meyn2012markov} are also not sufficient in this context.
However the stability of an alternative strategy where \eqref{eq:update-sigma} is replaced by a smooth update of the step-size
\del{\begin{equation}
	\label{eq:step-size-smooth}
	\sigma_{k+1} = \sigma_k \times \exp \left( \frac{1}{2 d_\sigma} \left( \frac{\mueff \|\sum_{i=1}^\mu w_i U_{k+1}^{s_{k+1}(i)}\|^2}{d}-1 \right) \right)
\end{equation}}%
has already been analyzed~\cite{toure2023global}.

As for CMA-ES, the following proposition gives a sufficient condition for assumption \Cref{ass:alpha} to hold. The proof goes as for \Cref{p:density-cma}, which can be found in \Cref{sec:additional-proofs}.
\begin{proposition}
\label{p:density-csa}
Suppose that $f$ has Lebesgue-negligible level sets. Define for all $z\in\cX$ and $v=(v_1,\dots,v_\mu)\in\cW$
\begin{equation}
	p_{z}(v)=\frac{\lambda!}{(\lambda-\mu)!} \1\left\{f\left(z+v_1\right)<\dots<f\left(z+v_\mu\right)\right\} \left(1-Q_{z}^f(v_\mu)\right)^{\lambda-\mu} \gamma^d(v_1)\dots \gamma^d(v_\mu)
\end{equation}
with $Q_{z}^f(u)=\int\1 \{f(z+\xi)<f(z+u)\} \gamma^d(\xi)\mathrm{d}\xi$ and where $\gamma^d$ is the density of the $d$-dimensionnal standard normal distribution w.r.t.\ Lebesgue. Then, $p_{z}$ defines a density (w.r.t.\ Lebesgue in $\bR^{d\mu}$) of the random variable $\alpha(z,U_1)$. Moreover, if $f$ is (a monotone transformation of) a continuous function, then $(z,v)\mapsto p_z(v)$ is l.s.c.
\end{proposition}

As for CMA-ES,
assumption \Cref{ass:F} holds since $F$, given in \eqref{eq:CSA-update-MC}, is the composition of a continuously differentiable function with the Lipschitz function $x\mapsto\|x\|$. Regarding \Cref{ass:sub_differe_S_steadily}, the next proposition states the existence of a steadily attracting state. The proof follows the same lines as \cite[Proposition 5.3]{chotard2019verifiable}, but is given for completeness.

\begin{proposition}
\label{p:steadily-attracting-csa}
Suppose that $f$ is continuous, scaling-invariant with Lebesgue-negligible level sets. Then, $0$ is a steadily attracting state.
\end{proposition}

\begin{proof}
For $z_0\in\bR^d$, we set $v_1=-[z_0,\dots,z_0]\in\bR^{d\mu}$, and $v_k=[0,\dots,0]\in\bR^{d\mu}$. Note that, by \Cref{p:density-csa}, since $f$ has Lebesgue-negligible level sets, $v_{1:k}\in\overline{\cO^k_{z_0}}$. Moreover, we have $S^k_{z_0}(v_{1:k})=0$ for every $k\geqslant1$, where $S^k_{z_0}$ is defined in \eqref{eq:ext-transition-map}. We conclude the proof by using \Cref{c:sufficient-condition-steadily-attracting}  .
\end{proof}

\del{For the controllabillity condition, as $F$ is locally Lipschitz but not continuously differentiable, we have to satisfy a more complex condition. However, to avoid to verify a full rank condition for all elements in $\partial S^1_{0} (v_1)$ for some $v_1\in\overline{\cO^{1}_{0}}$, }
To complete the verification of \Cref{ass:sub_differe_S_steadily}, we show in the next proposition that there exists $v_1\in\overline{\cO^{1}_{0}}$ such that $S^1_0$ is differentiable in $v_1$ and $\cD S^1_{0} (v_1)$ is of maximal rank.

\begin{proposition}
\label{p:rank-condition-csa}
Suppose that $f$ is continuous, scaling-invariant with Lebesgue-negligible level sets. Then, $S^1_{0}$ is differentiable in $v_1=(0,\dots,0)\in\overline{\cO^{1}_{0}}$ and $\cD S^1_{0}(v_1)$ is of maximal rank.
\end{proposition}

\begin{proof}
Note that, by \Cref{p:density-csa}, $v_1$ belongs to $\overline{\cO^{1}_{0}}$. Moreover, for $h=(h_1,\dots,h_\mu)\in\cW$, we have by definition of $F$ and of $S^1_0$, see \eqref{eq:CSA-update-MC} and \eqref{eq:ext-transition-map} respectively, that
\begin{align*}
	S^1_{0}(v_1+h) & = F(0,h) = \exp\left( \frac{1}{d_\sigma} \left( 1 - \frac{\sqrt{\mueff}\|\sum_{i=1}^\mu w_i h_i\|}{\bE\|\mathcal{N}(0,I_d)\|} \right) \right) \times \sum_{i=1}^\mu w_i h_i  .
\end{align*}
A simple Taylor expansion shows that
\begin{equation}
	\lim_{h\to 0} \frac{\left\| S^1_{0}(v_1+h) - S^1_{0}(v_1) -\exp\left(\frac{1}{d_\sigma}\right)\times \sum_{i=1}^\mu w_i h_i \right\|}{\|h\|}  = 0  ,
\end{equation}
ending the proof.
\end{proof}

Using \Cref{theo:1}, we deduce the $\varphi$-irreducibility and aperiodicity of the chain $\{z_k\}_{k\in\bN}$.
\begin{theorem}
Suppose that $f$ is continuous, scaling-invariant with Lebesgue-negligible level sets. Then, the Markov chain $\seq{z_k}{k\in\bN}$ defines a time-homogeneous $\varphi$-irreducible aperiodic T-chain, for which compact subsets of $\cX$ are small.
\end{theorem}
{
Note that in \cite{toure2023global}, it has been proven that the chain $\{z_k\}_{k\in\bN}$ is $\varphi$-irreducible, aperiodic and positive recurrent, on the condition that the step-size obeys to a smooth update instead of \eqref{eq:update-sigma}. However, a smooth step-size update was required only to prove the $\varphi$-irreducibility and aperiodicity of the chain, since 
the derivation of these two results rely in \cite{toure2023global} on results in \cite{chotard2019verifiable}. 
Now that we have proven that the chain $\{z_k\}_{k\in\bN}$ is $\varphi$-irreducible and aperiodic even when the step-size update is nonsmooth, we can prove that it is positive recurrent following
the proofs of
\cite{toure2023global}.}

\section{Proofs}\label{sec:proofs}~\!
\newcommand{\rankcondition}[1]{\textnormal{(}\hyperref[eq:full-rank-condition]{\ensuremath{\mathrm{R}_{#1}}}\textnormal{)}}%
{
We provide in this section the proofs of \Cref{theo:1} as well as of the results used to achieve the former.
While they are inspired from the previous works~\cite{meyn1991asymptotic,meyn2012markov,chotard2019verifiable},
the relaxation of the assumptions to state spaces being manifolds and an update function being locally Lipschitz represent a great challenge.
The manifold assumption requires to use at several places local arguments: an example is the proof of \Cref{l:small-set-lemma}, where we first prove the Euclidean case, and then we have to consider local charts to extend to manifolds.
The locally Lipschitz assumption brings other complications: as we cannot assume the differentiability at all states, we require a controllability condition---consisting in a full rank condition of all elements of the Clarke's derivative.
To this end,  we prove in \Cref{p:lsc-maxrank} the equivalence of the controllability condition \controllability{x} with a full rank condition \rankcondition{\bar x} at a neighbor point $\bar x$ where the update function is differentiable, based on Rademacher's theorem stating that a locally Lipschitz function possesses a dense set of points at which it is differentiable.
For the extension to carry over to manifolds, many tools from nonsmooth analysis need to be appropriately generalized.
Since this is not the main focus of this paper, it is relegated to \Cref{sec:appendix-clarke}. Finally, proofs that are straightforward adaptations of previous works are moved to \Cref{sec:additional-proofs}, where they are provided for completeness.
}




\newcommand{\CM}{\mathrm{CM}}

\subsection{Preliminary results}

\subsubsection{Accessibility, attracting and attainable states}

In this section, we generalize characterizations of globally {attracting states \cite{meyn1991asymptotic,meyn2012markov}} and steadily attracting states \cite{chotard2019verifiable}. 
In contrast to those previous references, we relax assumptions on the sets $\cX$, $\cU$ and $\cW$.
Indeed, \cite{meyn1991asymptotic,meyn2012markov,chotard2019verifiable} supposed that these sets were open subsets of Euclidean spaces.
Here, we only suppose that they are smooth connected manifolds, as formalized in \Cref{sec:model}. 
{
This generalization is relatively straightforward and as a result, their proofs  are given in Appendix~\ref{sec:additional-proofs} for completeness, as they are not the core of our contribution.}
For the rest of the paper, let us define $A_+^0(x)\coloneqq\{x\}$ and 
\begin{equation}
A_+^k(x) \coloneqq \{S_x^k(w_{1:k}) \mid w_{1:k}\in\cO_x^k\} \quad \text{for } k\geqslant1 .
\end{equation}
The set $A_+^k(x)$ is the set of states that can be reached by $\phi_k$ conditionally to $\phi_0=x$.
We also define the set of {attainable} states {\cite[Section~7.1.4]{meyn2012markov}, i.e.} that can be reached by $\seq{\phi_k}{k\in\bN}$ (in finite time) conditionally to $\phi_0=x$ as
\begin{equation}
\label{eq:A+-def}
A_+(x) \coloneqq \bigcup_{k\in\bN} A_+^k(x) .
\end{equation}
Then, we say that the control model associated to \eqref{eq:def-MC} is \textit{forward accessible} if for every $x\in\cX$, $A_+(x)$ has a nonempty interior in $\cX$ \cite{meyn2012markov}. 
Moreover, with these notations, a point $x^*\in\cX$ is a \textit{globally attracting state} \cite[Section~7.2.4]{meyn2012markov}\del{, see \Cref{sec:model}}, if for every $y\in\cX$ we have
\begin{equation}
\label{eq:def-globally-attracting-state}
x^*\in\bigcap_{T\geqslant1} \overline{\bigcup_{k\geqslant T}A_+^k(y)} .
\end{equation}
As shown in the next proposition which is exactly \cite[Proposition 3.1]{chotard2019verifiable} applied to our more general setting, this definition is equivalent to the statement we used in \Cref{sec:model} to introduce a globally attracting state that for any $y\in\cX$ and any neighborhood $U$ of $x^*$, there exists $k>0$ and a $k$-steps path between $y$ and $U$.
\begin{proposition}[Characterization of globally attracting states]
\label{p:characterization-globally-attracting}
Suppose \Cref{ass:alpha}. A point $x^*\in\cX$ is globally attracting if and only if one of the following equivalent conditions holds.
\begin{enumerate}
	\item[(i)] For any $y\in\cX$, $x^*\in\overline{A_+(y)}$.
	\item[(ii)] For any $y\in\cX$ and any open subset $U$ of $\cX$ containing $x^*$, there exist $k>0$ and a $k$-steps path from $y$ to $U$.
	\item[(iii)] For any $y\in\cX$, there exists a sequence $\{y_k\}_{k>0}$ with $y_k\in A_+^k(y)$, from which we can extract a subsequence converging to $x^*$.
\end{enumerate}
\end{proposition}
A point $x\in\cX$ is said to be \textit{reachable} by $P$ \cite[Section 6.1.2]{meyn2012markov} if for any measurable neighborhood $U$ of $x$ in $\cX$, we have
\begin{equation}
\forall y\in\cX, \quad \sum_{k\geqslant1} P^k (y,U)>0.
\end{equation}
The equivalence between globally attracting states and reachable states relies on the following proposition (see \cite[Proposition 3.2]{chotard2019verifiable}).
\begin{proposition}[Characterization of reachable states]
\label{p:equivalence-ksteps}
Consider the Markov kernel $P$ defined via \Cref{eq:def-MC}, and suppose \Cref{ass:alpha} and that $F$ is continuous. Then for any open subset $U$ of $\cX$, any $x\in\cX$ and $k>0$, the following statements are equivalent.
\begin{enumerate}
	\item[(i)] There exists a $k$-steps path from $x$ to $U$.
	\item[(ii)] $P^k(x,U)>0$.
\end{enumerate}
\end{proposition}
As an immediate consequence of \Cref{p:characterization-globally-attracting,p:equivalence-ksteps}, we get the following equivalence between states that are globally attracting by the control model associated to \eqref{eq:def-MC} and states that are reachable by $P$ (see \cite[Corollary~3.1]{chotard2019verifiable}).	
\begin{corollary}
\label{c:equivalence-attacting-reachable}
Consider the Markov kernel $P$ defined via \eqref{eq:def-MC}, and suppose \Cref{ass:alpha} and that $F$ is continuous. Then $x\in\cX$ is globally attracting if and only if it is reachable by $P$.
\end{corollary}

Recall that a state $x^*\in\cX$ is \textit{steadily attracting} \cite{chotard2019verifiable} if for all $y\in\cX$ and all open neighborhood $U$ of $x^*$ in $\cX$, there exists $T>0$ such that for all $k\geqslant T$ there exists a $k$-steps path from $y$ to $U$.

In the next proposition and corollary, we state two technical results related to steadily attracting states, which will be instrumental in the proofs of our main results. The next proposition is the equivalent for our setting of \cite[Proposition~3.3]{chotard2019verifiable}.
\begin{proposition}
\label{p:sufficient-condition-steadily-attracting}
Suppose \Cref{ass:alpha}. The following statements hold.
\begin{enumerate}
	\item[(i)] If $x^*\in\cX$ is steadily attracting, then it is globally attracting.
	\item[(ii)] A state $x^*\in\cX$ is steadily attracting if and only if for every $y\in\cX$ we can find a sequence $\{y_k\}_{k>0}$ with $y_k\in \overline{A_+^k(y)}$, which converges to $x^*$.
	\item[(iii)] Assume $F$ is continuous. If there exists a steadily attracting state, then every globally attracting state is steadily attracting.
\end{enumerate}
\end{proposition}
Note that the statement of \cite[Proposition~3.3~(ii)]{chotard2019verifiable} is slightly different as the element $y_k$ belongs to $A_+^k(y)$ while in (ii) above $y_k$ belongs to $\overline{A_+^k(y)}$. It is easy to see that both statements are equivalent.

In addition, we give the following corollary of \Cref{p:sufficient-condition-steadily-attracting} when $F$ is assumed continuous.
\begin{corollary}\label{c:sufficient-condition-steadily-attracting}
Suppose \Cref{ass:alpha} and that $F$ is continuous. Then, for any $x\in\cX$ and $k\in\bN$ we have the inclusion $\{S_x^k(w_{1:k}) \mid w_{1:k}\in \overline{\cO_x^k} \} \subset \overline{A_+^k(x)} $. Consequently the following statements are equivalent.
\begin{enumerate}
	\item[(i)] The point $x^*$ is steadily attracting.
	\item[(ii)] For every $x\in\cX$, there exists a sequence 
	$\{y_k\}_{k>0}$ satisfying $y_k\in\{S_x^k(w_{1:k}) \mid w_{1:k}\in \overline{\cO_x^k} \} $ for any $k\geqslant1$, and which converges to $x^*$.
	\item[(iii)] For every $x\in\cX$, for every neighborhood $U$ of $x^*$, there exists $T>0$ such that for any $k\geqslant T$ we can find $w_{1:k}\in\overline{\cO_x^k}$ satisfying $S_{x}^k(w_{1:k})\in U$.
\end{enumerate}
\end{corollary}

The next result corresponding to \cite[Proposition 3.4]{chotard2019verifiable} turns out to be  useful later in order to prove the aperiodicity of the Markov kernel $P$, given that it is $\varphi$-irreducible.
To this end, we need to introduce the notion of attainability, as considered in \cite{meyn2012markov}. We say that a state $x^*\in\cX$ is \textit{attainable} if 
\begin{equation}
\label{eq:attainable}
\forall y \in \cX,\quad x^*\in A_+(y).
\end{equation}
\begin{proposition}
\label{p:cycles}
Consider the Markov kernel $P$ defined via \eqref{eq:def-MC}, and suppose \Cref{ass:alpha}. Let $x^*\in\cX$ be attainable, and set
\begin{equation}
	E\coloneqq \{ a\in\bN^* \mid \exists T\in\bN, \forall k\geqslant T, x^*\in A_+^{ak}(x^*)\} .
\end{equation}
Then, the following statements hold.
\begin{enumerate}
	\item[(i)] $E$ is nonempty and for every $a, b\in E$, the greatest common divider of $a$ and $b$ satisfies $\gcd(a,b)\in E$.
	\item[(ii)] If $\gcd(E)=\max\{c\in\bN \mid c \text{ divides }a, \forall a \in E\}=1$, then $x^*$ is steadily attracting.
	\item[(iii)] If $P$ is $\varphi$-irreducible, then there exists a $d$-cycle {(as defined in \eqref{eq:aperiodicity-def})} with $d=\gcd(E)$.
\end{enumerate}
\end{proposition}

\subsubsection{Controllability condition}
In this section, we
{relax the Lipschitz assumption supposed in \cite{chotard2019verifiable} which was already a relaxation of the smooth assumption of \cite{meyn1991asymptotic,meyn2012markov}. {While many ideas are rooted in \cite{meyn1991asymptotic,meyn2012markov}, we} follow the exposition in \cite{chotard2019verifiable} and}
more precisely generalize \cite[Propositions 3.5, 3.6 and 3.7]{chotard2019verifiable}, 
{to obtain} condition \eqref{eq:controllabilty-condition}. 
The main challenge here is to deal with the condition that $F$ is supposed to be locally Lipschitz only.

First, as in \cite[Proposition~3.5]{chotard2019verifiable}, we prove that if the controllability condition \controllability{x^*} is satisfied for some $x^*$ a globally attracting state, then is satisfied for every $y\in\cX$.

\begin{proposition}
\label{p:rank-condition-on-all-states}
Suppose \Cref{ass:alpha} and \Cref{ass:F}. Let $x^*\in\cX$ be a globally attracting state. If \controllability{x^*} holds, then for every $y\in\cX$, \controllability{y} holds.
\end{proposition}

\begin{proof}
By \controllability{x^*}, 
there exist $k>0$ and $w^*_{1:k}\in\overline{\cO^k_{x^*}}$ such that $\partial S_{x^*}^k({w}^*_{1:k})$ is of rank $n$, the dimension of $\cX$.
See that, by \Cref{p:lsc-maxrank-supp},
we can assume that $w^*_{1:k}\in{\cO^k_{x^*}}$.
Moreover, the function $S^k:(z,w_{1:k})\mapsto S_z^k({w}_{1:k})$ is locally Lipschitz (since $F$ is locally Lipschitz), 
hence according to \Cref{p:usc-clarke-jacobian}, 
$\lim_{z \to x^*} \partial S^k(z,{w}_{1:k}) \subset \partial S^k(x^*,{w}_{1:k}) $ and since $\partial S^k(z,{w}_{1:k}) = \partial_z S^k(z,{w}_{1:k})\times\partial S_z^k({w}_{1:k})  $, we obtain $\lim_{z\to x^*} \partial S_z^k({w}^*_{1:k})\subset \partial S_{x^*}^k({w}^*_{1:k})$. 
Since $\rank$ is lower semicontinuous, 
we deduce that there exists an open neighborhood $U$ of $x^*$ such that for any $z\in U$,
$
\partial S_{z}^k({w}^*_{1:k})
$ 
is of rank $n$.
Moreover, $z\mapsto p_z^k({w}^*_{1:k})$ is lower semicontinuous, so, up to taking $U$ smaller, 
we can suppose that for any $z\in U$, $p_z^k({w}^*_{1:k})>0$, i.e., ${w}^*_{1:k}\in\cO_z^k$.

Let $y\in\cX$. Since $x^*$ is a globally attracting state, 
then by \Cref{p:characterization-globally-attracting}, 
there exist $t_0>0$ and $u_{1:t_0}$ a $t_0$-steps path from $y$ to $x\in U$, i.e., 
$u_{1:t_0}\in\cO^{t_0}_y$ and $x= S_y^{t_0}(u_{1:t_0})\in U$. 
Since  $\cO^{t_0}_y$ is open and $S_y^{t_0}$ is continuously locally Lipschitz, 
by \Cref{c:rademacher}
we can assume w.l.o.g.\ that $S_y^{t_0}$ is differentiable at $u_{1:t_0}$.

Since $x\in U$, then $\partial S_{x}^k({w}^*_{1:t_0})$ is of maximal rank, using the chain rule, see \Cref{p:chain-rule-manifolds}, we deduce that, for $T=t_0+k$ and $u_{t_0+1:t_0+k}=w^*_{1:k}$, we have that $\partial S_y^T(u_{1:T})$ is of maximal rank.
\end{proof}

The next proposition states that if we find a point $x^*$, $k>0$ and ${w}_{1:k}^*\in\cO_{x^*}^k$ which satisfy the forementionned controllability condition \controllability{x^*}, that is, $\partial S_{x^*}^k ({w}_{1:k}^*)$ is of maximal rank, then, using \Cref{p:lsc-maxrank}, we can find ${u}_{1:k}^*\in\cO_{x^*}^k$ as closed as we want from ${w}_{1:k}^*$ such that $S_{x^*}^k$ is differentiable in ${u}_{1:k}^*$ and $\cD S_{x^*}^k({u}_{1:k}^*)$ is of maximal rank. In other words, our controllability condition \controllability{x^*} implies a full rank condition. 

\begin{proposition}
\label{p:steadily-att-state-differentiable}
Suppose \Cref{ass:alpha} and \Cref{ass:F}. Let $x^*\in\cX$ and suppose that \controllability{x^*} holds. Then, condition \rankcondition{x^*} stated below holds.
\end{proposition}

\begin{proof}
By \controllability{x^*} and by \Cref{p:lsc-maxrank-supp}, there exist $k>0$ and ${w}^*_{1:k}\in\cO_{x^*}^k$ such that $\partial S_{x^*}^k({w}^*_{1:k})$ is of maximal rank. 
By \Cref{p:lsc-maxrank} below, for any neighborhood $W\subset\cW$ of ${w}^*_{1:k}$, there exists ${u}_{1:k}^*\in W$, such that $S_{x^*}^k$ is differentiable in ${u}_{1:k}^*$, with $\rank\cD S_{x^*}^k({u}_{1:k}^*)=n$. However $\cO_{x^*}^k$ is open, so we can take $W=\cO_{x^*}^k$ and complete the proof.
\end{proof}
\begin{proposition}
\label{p:lsc-maxrank}
Suppose that $f\colon \cX\to\cY$ is locally Lipschitz at $x_0\in\cX$, and that $\partial f (x_0)$ is of maximal rank, i.e., any $h\in\partial f (x_0)$ is of maximal rank. Then, there exists a neighborhood $U$ of $x_0$ such that for any $y\in U$, $\partial f(y)$ is of maximal rank. Moreover, for every neighborhood $V\subset U$ of $x_0$, there exists $y_0\in V$ such that $f$ is differentiable at $y_0$ and $\cD f(y_0)$ is of maximal rank.
\end{proposition}

\begin{proof}
Let $A=\{h\in\cL(T_{x_0}\cX,T_{f(x_0)}\cY)\mid h \text{ is not of maximal rank}\}$. Since the application $\rank$ is l.s.c., then $A$ is a closed set. By \Cref{p:compact-clarke-gradient}, $\partial f(x_0)$ is compact, and disjoint from $A$ since it is assumed to be of maximal rank. Thus $\mathrm{dist}(\partial f(x_0), A)>0$, where $\mathrm{dist}$ is a metric induced by a norm on the finitely dimensioned affine space $\cL(T_{x_0}\cX,T_{f(x_0)}\cY)$. Moreover, there exists $h^*\in\partial f(x_0)$ such that for every $h\in\partial f(x_0)$ we have
$$
\mathrm{dist}(h,A)\geqslant \mathrm{dist}(h^*,A)=\mathrm{dist}(\partial f(x_0), A) >0.
$$
By \cite[Proposition 2.6.2(c)]{clarke1990optimization}, there exists a neighborhood $U$ of $x_0$ such that for all $y\in U$, $\mathrm{dist}(\partial f(y),A) $ $\geqslant \mathrm{dist}(h^*,A)/2 >0$, thus $\partial f(y)$ is of maximal rank. The second part follows from Rademacher's theorem, see \Cref{c:rademacher}.
\end{proof}
From now on, we can assume a full rank condition, i.e.,
\begin{equation}
\tag{\ensuremath{\mathrm{R}_x}}
\label{eq:full-rank-condition}
\text{there exists } w_{1:k}\in{\cO^{k}_{x}} \text{ such that } \cD S_x^k (w_{1:k}) \text{ exists and is of maximal rank,}
\end{equation}
instead of the controllability condition \eqref{eq:controllabilty-condition}. We can then use \Cref{p:steadily-att-state-differentiable} to extend our results. {The next proposition states that if we can find a globally attracting state $x^*$ satisfying the maximal rank condition \rankcondition{x^*}, then we can find an attainable state.} {It generalizes \cite[Proposition~3.6]{chotard2019verifiable}.}

\begin{proposition}
\label{p:implicit-Skx}
Suppose \Cref{ass:alpha} and \Cref{ass:F}. Let $x^*\in\cX$ and suppose that there exist $k>0$ and ${w}_{1:k}^*\in\cO_{x^*}^k$ such that \rankcondition{x^*} is satisfied with ${w}_{1:k}^*$. 
\begin{enumerate}
	\item[(i)] There exists $U$ a neighborhood of $x^*$ such that for any $x\in U$, there exists ${w}_{1:k}\in\cO_x^k$ for which $S_x^k({w}_{1:k})=S_{x^*}^k({w}_{1:k}^*)$.
	\item[(ii)] If $x^*$ is globally attracting, then $S_{x^*}^k({w}_{1:k}^*)$ is attainable{, see \eqref{eq:attainable}}.
\end{enumerate}
\end{proposition}

\begin{proof}
(i) Let $(U,\varphi)$ be a local chart of $\cX$ around $x^*$, $(V,\theta)$ a local chart of $\cX$ around $S_{x^*}^k({w}_{1:k}^*)$, and $(W,\psi)$ a local chart of $\cW^k$ around ${w}_{1:k}^*$, such that the following differentiable function is well-defined
$$
\tilde{S}^k\colon \begin{array}{l}
	\varphi(U)\times\psi(W)\subset\bR^n \times \bR^{kp} \to V\subset\bR^n \\
	(x,w)=((x_1,\dots,x_n),(w_1,\dots,w_{kp})) \mapsto \tilde{S}^k_{x} (w) \coloneqq \theta\circ S^k_{\varphi^{-1}(x)}\circ\psi^{-1}(w).
\end{array}
$$
We recall that the positive integers $n$ and $p$ are the dimensions of $\cX$ and $\cW$, respectively.
By composition, observe that $\cD_w \tilde{S}^k (\varphi(x^*),\psi({w}_{1:k}^*))$ is surjective. Hence, we can find coordinates $i_1,\dots,i_n$ of $\bR^{kp}$ such that 
$$
\det\left[\cD_{w_{i_1}} \tilde{S}^k \mid \dots \mid \cD_{w_{i_n}}\tilde S^k \right] (\varphi(x^*),\psi({w}_{1:k}^*)) = n .
$$
Note that, up to a permutation of indices in the chart $\psi$, we can assume w.l.o.g.\ that $i_1,\dots,i_n$ equal respectively $kp-n+1,\dots,kp$. 
To ease the presentation,
we use the following abuse of notation $(w_1^*,\dots,w_{kp}^*)=\psi({w}_{1:k}^*)$. Then, by the implicit function theorem, see \Cref{t:implicit-function-theorem}, there exist neighborhoods $M$ of $(\varphi(x^*),w^*_1,\dots,w^*_{kp-n})$ and $N$ of $(w^*_{kp-n+1},\dots,w^*_{kp})$, and a $\cC^1$ function $g\colon M\to N$ such that, for every $(x_1,\dots,x_n,$ $w_{1},\dots,w_{kp-n})\in M$, we have
$$
\tilde{S}^k_{(x_1,\dots,x_n)} \left(w_1,\dots,w_{kp-n},g(x_1,\dots,x_n,w_{1},\dots,w_{kp-n})\right) = \tilde{S}^k_{\varphi(x^*)}(w_1^*,\dots,w_{kp}^*) .
$$
This proves (i).

{
	(ii) Suppose that $x^*$ is globally attracting.
	Let $U\subset\cX$ be a neighborhood of $x^*$ satisfying (i), and let $y\in\cX$.
	Then, by \Cref{p:characterization-globally-attracting}(ii), there exist $k_1>0$ and $w_{1:k_1}\in\cO_y^{k_1}$ such that $S_y^{k_1}(w_{1:k_1})\in U$. 
	Since $U$ satisfies (i), there exists $w_{k_1+1:k_1+k}\in\cO_{S_y^{k_1}(w_{1:k_1})}^k$ with $S_y^{k_1+k}(w_{1:k_1+k})=S_{x^*}^k(w_{1:k}^*)$.
}
\end{proof}

We discuss in the next proposition generalizing  \cite[Proposition~3.7]{chotard2019verifiable} the forward accessibility of the control model \eqref{eq:def-MC}. 
We recall that it is said to be forward accessible if for every $x\in\cX$, the subset $A_+(x)\subset\cX$ defined in \eqref{eq:A+-def} of states that can be reached in finite time starting from $x$, has a nonempty interior.

\begin{proposition}
\label{p:forward-accessibility}
Suppose \Cref{ass:alpha} and \Cref{ass:F}. If for every $x\in\cX$, \rankcondition{x} holds, then the control model associated to \eqref{eq:def-MC} is forward accessible.

Furthermore, if $F$ is smooth (infinitely differentiable), the control model is forward accessible if and only if for every $x\in\cX$, \rankcondition{x} holds.
\end{proposition}

\begin{proof}
We apply the Local Submersion Theorem~\cite[Chapter 1.4]{guillemin2010differential}. Since $S_{x}^k$ is a submersion at $w_{1:k}$, there exist local charts $(W,\psi)$ of $\cW^k$ around $w_{1:k}$ and $(V,\varphi)$ of $\cX$ around $S_{x}^k(w_{1:k})$ such that 
$$
\varphi\circ S_{x}^k \circ \psi (u_1,\dots,u_{kp}) = (u_1,\dots,u_n) \quad \text{for all } (u_1,\dots,u_{kp})\in\psi(W).
$$
Therefore, since $\varphi$ is a continuous bijection (by definition of a local chart), then there exists a neighborhood $U$ of $S_{x}^k(w_{1:k})$ such that $S_{x}^k(W)=U$. Moreover, $\cO_{x^*}^k$ is an open subset of $\cW^k$, so we can assume $W\subset\cO_{x}^k$. Therefore, $U\subset A_+^k(x)$, which hence has a nonempty interior.

Suppose now that $F$ is smooth and that the control model is forward accessible. Then, for every $x\in\cX$, $\mathrm{int}(A_+(x))\neq\emptyset$.
Since $A_+(x)=\cup_{k\geqslant0}A_+^k(x)$, we deduce that there exists $k\in\bN$ such that $\mathrm{int}(A_+^k(x))\neq\emptyset$. 
Since $\mathrm{int}(A_+^0(x))=\mathrm{int}(\{x\})=\emptyset$, we find that necessarily $k>0$. By Sard's theorem~\cite[Appendix 1]{guillemin2010differential}, we have that the set $N\coloneqq\{\mathbf{w}\in\cO_x^k \mid \rank\cD S_x^k(\mathbf{w})<n\}$ is of measure zero, that is, for all charts $(\varphi,U)$ of $\cX$, we have $\Leb~ \varphi(N\cap U)=0$, hence $\mathrm{int}(N)=\emptyset$. We deduce that there exists $w_{1:k}\in\cO_x^k\setminus N$, i.e., such that $\rank\cD S_x^k(w_{1:k})=n$.
\end{proof}

\subsection{Proofs of the main results: verifiable conditions for irreducibility and aperiodicity}\label{sec:proofs-main-results}

\subsubsection{T-chain and irreducibility}
We preface our proofs by an extension of  \cite[Lemma~3.0]{meyn1991asymptotic} to our context, that is, for a locally Lipschitz function between manifolds instead of a smooth function between open subsets of Euclidean spaces. 

\begin{lemma}
\label{l:small-set-lemma}
Let $\cX_1$ be a $n$-dimensional manifold, $\tilde{\cW}_1$ a $m$-dimensional manifold, $\hat{\cW}_1$ a $n$-dimensional manifold, equipped with their respective Borelian $\sigma$-fields and with a measure $\zeta_{\cX}$ (resp.\ $\zeta_{\tilde{\cW}}$, $\zeta_{\hat{\cW}}$), which satisfies that for any $A\in\cB(\cX_1)$ (resp.\ of $\cB(\tilde{\cW}_1)$, $\cB(\hat{\cW}_1)$), $\zeta_{\cX}(A)=0$ (resp.\ $\zeta_{\tilde{\cW}}(A)=0$, $\zeta_{\hat{\cW}}(A)=0$) if and only if $\varphi(A\cap U)$ is Lebesgue-negligible for every chart $(\varphi,U)$. 

Let $G\colon(x,\tilde w,\hat w)\in\cX_1\times\tilde{\cW}_1\times\hat{\cW}_1\mapsto z\in\cX_1$ be a locally Lipschitz map differentiable in $(x_0,\tilde w_0,\hat w_0)$ such that $\rank\cD_{\hat{w}}G(x_0,\tilde{w}_0,\hat{w}_0)=n$. Then,
\begin{enumerate}
	\item[(i)] There exists an open subset $\cX\times\tilde{\cW}\times\hat{\cW}\subset\cX_1\times\tilde{\cW}_1\times\hat{\cW}_1$ containing $(x_0,\tilde{w}_0,\hat{w}_0)$ such that for any $x\in\cX$, the measure defined by
	\begin{equation}
		\nu(x,\cdot)\colon A\subset \cX_1 \mapsto \int_{\tilde\cW}\int_{\hat{\cW}} \1_A \{G(x,\tilde{w},\hat{w})\} d\zeta_{\tilde{\cW}}(\tilde{w})d\zeta_{\hat{\cW}}(\hat{w})
	\end{equation}
	is equivalent to the measure $\zeta_{\cX}$ on an open subset $\cR_x$ of $\cX_1$.
	\item[(ii)] There exist $c>0$, $U_{x_0}$ an open subset of $\cX_1$ containing $x_0$, $V_{x_0}^{\tilde{w}_0,\hat{w}_0}$ an open subset of $\tilde{\cW}_1$ containing $G(x_0,\tilde{w}_0,\hat{w}_0)$ such that for every $x\in\cX$ and every measurable subset $A$ of $\cX_1$, we have $\nu(x,A)\geqslant c\1_{U_{x_0}}(x)\times\zeta_{\cX_1}(A\cap V_{x_0}^{\tilde{w}_0,\hat{w}_0})$.
\end{enumerate}
\end{lemma}

{The proof of this result proceeds in two steps. First, we assume the spaces to be Euclidean (as in \cite{meyn1991asymptotic}) while allowing the function to be only locally Lipschitz. This requires the use of more general results applicable to nonsmooth functions~\cite{clarke1990optimization}, as well as a general change-of-variable property~\cite[Theorem 3]{hajlasz1993change}. The second step extends the proof to manifolds.
}

\begin{proof}
First we prove the lemma when $\cX_1,\tilde{\cW}_1,\hat{\cW}_1$ are open subsets respectively of $\bR^n$, $\bR^m$, $\bR^n$, and $\zeta_{\cX},\zeta_{\tilde{\cW}},\zeta_{\hat{\cW}}$ are assumed to be the Lebesgue measures on $\bR^n$, $\bR^m$, $\bR^n$ respectively. 

Define the function
$$
G^\star \colon(x,\tilde w,\hat w)\in\cX_1\times\tilde{\cW}_1\times\hat{\cW}_1\mapsto (x,\tilde w , G(x,\tilde w,\hat w)) \in \bR^{n+m+n}
$$
Then, since $\cD_{\hat{w}}G(x_0,\tilde w_0,\hat w_0)$ is of rank $n$, then $\cD G^\star (x_0,\tilde w_0,\hat w_0)$ exists and is a full-rank squared matrix. Therefore, the inverse function theorem --as stated in \Cref{t:inverse-function-theorem}-- applies and we find a neighborhood $\cX\times\tilde\cW\times\hat\cW$ of $(x_0,\tilde w_0,\hat w_0)$, a neighborhood $\cR$ of $(x_0,\tilde w_0,z_0)$ (where $z_0\coloneqq G(x_0,\tilde w_0,\hat w_0)$), and a locally Lipschitz function $H^\star \colon\cR\to\cX\times\tilde\cW\times\hat\cW$ such that
$$
H^\star(G^\star(x,\tilde w,\hat w)) = (x,\tilde w,\hat w) \quad \text{for every } (x,\tilde w,\hat w) \in \cX\times\tilde\cW\times\hat\cW .
$$
Thus, there exists a locally Lipschitz function $H\colon\cR\to\hat{\cW}$ such that
$$
H(x,\tilde w,G(x,\tilde w,\hat w)) = \hat w \quad \text{for every } (x,\tilde w,\hat w) \in \cX\times\tilde\cW\times\hat\cW .
$$
Then, by the chain rule, see \cite[Theorem 2.6.6]{clarke1990optimization}, for every $(x,\tilde w,\hat w) \in \cX\times\tilde\cW\times\hat\cW$ at which $G$ admits a partial derivative w.r.t.\ $\hat w$, we have that
\begin{equation}
	\label{eq:inverse-jacobian-H}
	\cD_z H (x,\tilde w, G(x,\tilde w,\hat w)) = \left[ \cD_{\hat{w}} G(x,\tilde w,\hat w) \right]^{-1}
\end{equation}
which is thus invertible. Moreover, by \cite[Proposition 2.6.2(c)]{clarke1990optimization}, $\cD_z H$ is continuous at points on which it is defined (which is dense by Rademacher's theorem~\cite[Theorem 3.2]{gariepy2015measure}). Therefore, there exists $h_0>0$ such that in each of these points, by \eqref{eq:inverse-jacobian-H} we have
\begin{equation}
	\label{eq:minoration-detDH}
	\left| \det \cD_z H \right| \geqslant h_0 .
\end{equation}
Then, applying \cite[Theorem 3]{hajlasz1993change} and Fubini's theorem, we get
\begin{equation}
	\nu(x,A) = \int_{A} \left( \int \1_{\cR} (x,\tilde w,z)\left| \det \cD_z H \right| \,\d\tilde w \right) \,\d z , 
\end{equation}
so that 
\begin{equation}
	p(x,z) \coloneqq \int \1_{\cR} (x,\tilde w,z)\left| \det \cD_z H \right| \,\d\tilde w 
\end{equation}
defines a density w.r.t.\ Lebesgue for $\nu(x,\cdot)$. The rest of proof goes as in \cite[Lemma 3.0]{meyn1991asymptotic}, that we recall here for completeness.

Fix $x\in\cX$ and let $\cR_x$ be the open subset of $\bR^n$ defined by
$$
\cR_x= \left\{ z\in\cX_1 \mid \exists \tilde w \in\tilde\cW , (x,\tilde w,z)\in\cR \right\} .
$$ 
Then, note that $p(x,z)$ is positive if and only if $z\in\cR_x$, and zero otherwise. This proves (i). For (ii), observe that, since $\cR$ is a neighborhood of $(x_0,\tilde w_0,z_0)$, then it contains a nonempty open subset $\cX_0\times\tilde\cW_0\times\msz_0$ containing $(x_0,\tilde w_0,z_0)$. We get then that $p(x,z)\geqslant h_0 \times \Leb(\tilde\cW_0)$ for every $(x,z)\in\cX_0\times\cZ_0$. Then,
$$
\nu(x,A) \geqslant h_0\Leb(\tilde\cW_0)\1\{x\in\cX_0\}\times\Leb(A\cap\cZ_0)
$$
which proves (ii).

Now suppose that $\cX_1,\tilde{\cW}_1,\hat{\cW}_1$ are manifolds.

Let $(\varphi,\cX_2)$ be a local chart of $\cX_1$ around $x_0$, $(\tilde \psi,\tilde\cW_2)$ be a local chart of $\tilde\cW_1$ around $\tilde w_0$, $(\hat \psi,\hat\cW_2)$ be a local chart of $\hat\cW_1$ around $\hat w_0$, and $(\eta,\cX_3)$ be a local chart of $\cX_1$ around $z_0=G(x_0,\tilde w_0,\hat w_0)$.

Then, define the locally Lipschitz map
$$
G^{\mathrm{loc}}\colon (x,\tilde w,\hat w) \in \varphi(\cX_2)\times\tilde{\psi}(\tilde\cW_2)\times\hat{\psi}(\hat\cW_2)\mapsto z = \eta\circ G(\varphi^{-1}(x),\tilde{\psi}^{-1}(\tilde w),\hat{\psi}(\hat w)) \in \bR^n .
$$
Thus, (i) and (ii) hold with $G^{\mathrm{loc}}$, and 
$$
\nu^{\mathrm{loc}} (x,A) = \int_{\tilde \cW_0} \int_{\hat \cW_0} \1_A (x,\tilde w,\hat w) \eta\circ G(\varphi^{-1}(x),\tilde\psi^{-1}(\tilde w),\hat\psi^{-1}(\hat w)) \, \d\tilde w \,\d\hat w
$$
is equivalent to the Lebesgue measure, for all $x\in\cX_0$, and $\cX_0\times\tilde\cW_0\times\hat\cW_0$ being a neighborhood of $(\varphi(x_0),\tilde\psi(\tilde w_0),\hat\psi(\hat w_0))$. But, by assumption on the measures $\zeta_{\cX}$, $\zeta_{\tilde{\cW}}$ and $\zeta_{\hat{\cW}}$, $\nu(x,\cdot)$ is locally equivalent to $\rho^{-1}\circ\nu^{\mathrm{loc}} (\varphi(x),\cdot)$ for all local chart $(\rho,A)$ of $\cX_1$, thus is locally equivalent to $\eta^{-1}\circ\Leb_n$ where $\Leb_n$ is the Lebesgue measure of $\bR^n$. Thus, $\nu(x,\cdot)$ is equivalent to $\zeta_{\cX}$. This proves (i).

Now apply (ii) to $G^{\mathrm{loc}}$, and find $c>0$, $U_{x_0}$ an open of $\varphi(\cX_2)$ containing $\varphi(x_0)$, $V_{x_0}^{\tilde w_0,\hat w_0}$ an open of $\cX_3$ containing $\eta(z_0)$, such that
$$
\nu^{\mathrm{loc}}(x,A) \geqslant c\1_{U_{x_0}}(x)\times\Leb_n(A\cap V_{x_0}^{\tilde w_0,\hat w_0}) \quad \text{for every } x\in\varphi(\cX), A\subset \bR^n.
$$
But, by assumption on $\zeta_{\cX}$, we find $L_1^\eta,L_2^\eta>0$ such that
\begin{align*}
	\nu (x,A) & \geqslant L_1^\eta\times\nu^{\mathrm{loc}}(\varphi(x),\eta(A\cap\cX_3)) \\
	& \geqslant L_1^\eta \times c\1_{U_{x_0}}(\varphi(x))\times \Leb_n(\eta(A\cap\cX_3)\cap V_{x_0}^{\tilde w_0,\hat w_0}) \\
	& \geqslant L_1^\eta \times c\1_{\varphi^{-1}(U_{x_0})}(x) \times L_2^\eta \times \zeta_{\cX} (A\cap \eta^{-1}(V_{x_0}^{\tilde w_0,\hat w_0}))
\end{align*}
for all $x\in\cX$ and $A\subset \cX_1$, which proves (ii).
\end{proof}

We can now state the following result.

\begin{proposition}
\label{p:smallset-property}
Consider the Markov kernel $P$ defined via \eqref{eq:def-MC}, and suppose \Cref{ass:alpha}-\Cref{ass:F}. Let $x\in\cX$. 
\begin{enumerate}
	\item[(i)] If \rankcondition{x} holds for some $k>0$ and ${w}_{1:k}\in\cO_x^k$, then there exist $c>0$, and open subsets $U_x$ and $V_x^{{w}_{1:k}}$ of $\cX$ containing $x$ and $S_x^k({w}_{1:k})$ respectively, such that
	\begin{equation}
		P^k(y,A)\geqslant c \zeta_{\cX} \left( A \right) \quad \text{for every } y\in U_x \text{ and } A \in \cB(\cX) ,
	\end{equation}
	for some nontrivial measure $\zeta_{\cX}$ on $V_x^{{w}_{1:k}}$.
	That is, $U_x$ is a $k$-small set.
	\item[(ii)] If furthermore $F$ is smooth (infinitely differentiable), and if there exist $k>0$, $c>0$ and $(\varphi,V)$ a local chart of $\cX$ such that
	\begin{equation}
		\label{eq:small-set2}
		P^k(x,A) \geqslant c\Leb\circ\varphi(A\cap V) \quad \text{for every } A \in\cB(\cX) ,
	\end{equation}
	then \rankcondition{x} holds.
\end{enumerate}
\end{proposition}

\begin{proof}
Condition \rankcondition{x} implies that
$\rank \cD S_{x}^k(w_{1:k})=n$ for some $k>0$ and $w_{1:k}\in\cO_{x}^k$. Since $(\bar x,\bar w_{1:k})\mapsto p_{\bar x} ^k(\bar w_{1:k})$ is l.s.c., and $p_{x}^k(w_{1:k})>0$, then there exist $p_0>0$ and a neighborhood $\cX_1\times\cW_1$ of $(x,w_{1:k})$ such that $p_{\bar x}^k(\bar w_{1:k})\geqslant p_0$ for every $\bar x\in\cX_1$ and every $\bar w_{1:k}\in\cW_1$. Then, for every $y\in\cX_1$, we have
\begin{align}
	\label{eq:minoration-kernel}
	P^k(y,A)  = \int_{\cO_y^k} \1_A(S_y^k(\bar w_{1:k})) p_y^k (\bar w_{1:k}) \d\zeta_{\cW}^{\otimes k}(\bar w_{1:k}) \geqslant p_0 \int_{\cW_1}\1_A(S_y^k(\bar w_{1:k}))\d\zeta_{\cW}^{\otimes k}(\bar w_{1:k}) .
\end{align}
Since $S_{x}^k$ is a submersion at $w_{1:k}$, by the Local Submersion theorem, 
there exists a local chart $(V,\psi)$ of $\cW^k$ around $w_{1:k}$ and a local chart $(U,\varphi)$ of $\cX$ around $S_{x}^k(w_{1:k})$, such that for every $(\bar w_1,\dots,\bar w_{kp})\in\psi(V)$, we have 
\begin{equation}
	\label{eq:local-sub-th-Skx}
	\varphi\circ S_{x}^k\circ \psi^{-1} (\bar w_1,\dots,\bar w_{kp}) = (\bar w_1,\dots,\bar w_n).
\end{equation}
Note that, up to taking $V$ and $\cW_1$ smaller, we can assume $V=\cW_1${,
	and that $\psi(\cW_1)=\hat \cW_1\times\tilde \cW_1$ is a rectangle of $\bR^{kp}$,
	with $\hat \cW_1=\{(\bar w_1,\dots,\bar w_n)\in\bR^n \mid \exists (\bar w_{n+1},\dots,\bar w_{kp})\in\bR^{kp-n}, (\bar w_1,\dots,\bar w_{kp})\in\psi(\cW_1)\}\subset \bR^n$ and likewise $\tilde\cW_1\subset\bR^{kp-n}$.
	Hence, the function 
	\begin{equation}
		\begin{array}{rl}
			G\colon~&\cX\times\hat\cW_1\times\tilde\cW_1 \to \cX \\
			&(\bar x,(\bar w_1,\dots,\bar w_n),(\bar w_{n+1},\dots,\bar w_{kp})) \mapsto S_{\bar x}^k \circ \psi^{-1} (\bar w_1,\dots,\bar w_{kp})
		\end{array}
	\end{equation} 
	satisfies, by \eqref{eq:local-sub-th-Skx}, that $\rank \cD_{(w_1,\dots,w_n)}G(x,w_1,\dots,w_{kp})=n$.} 
Then, by \Cref{l:small-set-lemma} and \Cref{ass:alpha}(ii), there exist $c>0$, $U_{x}$ an open subset of $\cX$ containing $x$, $(\varphi,V_{x}^{w_{1:k}})$ a local chart of $\cX$ around $S_{x}^k(w_{1:k})$ such that, for every $y\in\cX$, we have
\begin{equation}
	\label{eq:minoration-lebesgue}
	\int_{\cW_1} \1_A(S_y^k(\bar w_{1:k})) \d\zeta_{\cW}^k(\bar w_{1:k}) \geqslant c \1_{U_{x}}(y) \zeta_{\cX} (A\cap V_{x}^{w_{1:k}}) ,
\end{equation}
with $\zeta_{\cX} = \Leb\circ\varphi(\cdot \cap V_{x}^{w_{1:k}})$ is a measure which satisfies the assumption required in \Cref{l:small-set-lemma} on $V_{x}^{w_{1:k}}$.
Combining \eqref{eq:minoration-kernel} and \eqref{eq:minoration-lebesgue} gives $P^k(y,A) \geqslant c p_0 \zeta_{\cX} (A\cap V_{x}^{w_{1:k}})$ for every $y\in U_{x}\cap \cX_1$, proving (i).

Suppose now that $F$ is smooth, then $(\bar x,\bar w_{1:k})\mapsto S_{\bar x}^k(\bar w_{1:k})$ is smooth for all $\t>0$.
Take $\t>0$, $c>0$ and $(\varphi,V)$ a local chart such that \eqref{eq:small-set2} holds.
Let $N=\{S_x^k(\bar w_{1:k})\in\cX\mid \bar w_{1:k}\in\cO_x^k, \rank\cD S_x^k(\bar w_{1:k})<n\}$. By Sard's theorem, we know that $\Leb\circ\varphi(N\cap V)=0$, implying that $P^k(x,V\setminus N)\geqslant c\Leb\circ\varphi(V\setminus N)= c\Leb\circ\varphi(V)>0$. Hence there exists $w_{1:k}\in\cO_{x}^k$ such that $S_x^k(w_{1:k})\in V\setminus N$, i.e., $\rank\cD S_x^k(w_{1:k})=n$.
\end{proof}

Following \cite[Corollary 4.1]{chotard2019verifiable}, we now deduce sufficient conditions for the Markov kernel $P$ to define a T-chain.

\begin{corollary}
\label{c:Tchain}
Consider the Markov kernel $P$ defined via \eqref{eq:def-MC}, and suppose \Cref{ass:alpha} and \Cref{ass:F}. Suppose that for any $x\in\cX$, \controllability{x} holds. Then $\cX$ can be written as the union of open small sets and thus $P$ is a T-chain.
\end{corollary}

\begin{proof}
First, using \Cref{p:steadily-att-state-differentiable}, for all $x\in\cX$, \rankcondition{x} holds, i.e., there exist $k>0$ and ${u}_{1:k}\in\cO_{x}^k$ such that $S_x^k$ is differentiable in ${u}_{1:k}$ and $\rank \cD S_x^k({u}_{1:k})=n$.

\Cref{p:smallset-property} implies that for every $x\in\cX$, there exists an open neighborhood $\cU_x$ of $x$ in $\cX$ which is a $k$-small set. Denoting $a$ the Dirac distribution in $k$, we find that $\cU_x$ is $\nu_a$-petite, hence, by \cite[Proposition 6.2.3]{meyn2012markov}, $K_a$ possesses a continuous component $T$ which is nontrivial on $\cU_x$ and in particular at $x$. Thus, by \cite[Proposition 6.2.4]{meyn2012markov}, $P$ is a T-chain.
\end{proof}

We now characterize the support of the maximal irreducibility measure of $P$. 
{We recall that, by \cite[Proposition 4.2.2]{meyn2012markov}, any $\varphi$-irreducible Markov kernel $P$ admits a maximal irreducibility measure $\psi$, that is, $P$ is $\psi$-irreducible and for every irreducibility measure $\varphi$ of $P$, we have that $\mathrm{supp}~\varphi\subset\mathrm{supp}~\psi$.}
The proof mimics the one of \cite[Proposition 4.2]{chotard2019verifiable},
and is given for completeness in \Cref{sec:additional-proofs}.
\begin{proposition}
\label{p:support-irreducibility-measure}
Suppose that $P$ is a $\psi$-irreducible Markov kernel, defined via \eqref{eq:def-MC}, with $\psi$ a maximal irreducibility measure, that \Cref{ass:alpha} holds and that $F$ is continuous. Then
\begin{equation}
	\label{eq:supp-psi=globally-attracting-states}
	\mathrm{supp}~\psi=\{x^*\in\cX \mid x^* \text{ is globally attracting}\} .
\end{equation}
Furthermore, if $x^*\in\cX$ is globally attracting, then
\begin{equation}
	\label{eq:supp-psi=A+}
	\mathrm{supp}~\psi=\overline{A_+(x^*)}.
\end{equation}
\end{proposition}

We now state our core results, 
from which we deduce \Cref{theo:1}. 
Assuming the controllability condition is satisfied at every $x$, there is equivalence between the irreducibility of $P$ and the existence of a globally attracting state.
\begin{theorem}
\label{t:sufficient-condition-ireducibility}
Consider the Markov kernel $P$ defined via \eqref{eq:def-MC}, and suppose \Cref{ass:alpha} and \Cref{ass:F}. Suppose \controllability{x} is satisfied for every $x\in\cX$. Then $P$ is $\varphi$-irreducible if and only if a globally attracting state exists.
\end{theorem}
\begin{proof}
By \Cref{p:steadily-att-state-differentiable}, we know that \rankcondition{x} holds for any $x\in\cX$. 
If $P$ is $\varphi$-irreducible, then by \Cref{p:support-irreducibility-measure}, any point of the support of the nontrivial measure $\varphi$ is globally attracting, hence there exists a globally attracting state. Conversely, if $x^*\in\cX$ is globally attracting, then, by \Cref{c:equivalence-attacting-reachable}, $x^*$ is reachable by $P$,
and by \Cref{c:Tchain}, $P$ is a T-chain.
As a result, by \cite[Proposition 6.2.1]{meyn2012markov}, $P$ is $\varphi$-irreducible.
\end{proof}

We deduce from this theorem our first practical result in order to prove the irreducibility, the T-chain property of a Markov kernel following the model investigated. 
If assumptions \Cref{ass:alpha} and \Cref{ass:F} are satisfied for a Markov kernel defined via \eqref{eq:def-MC}, the theorem below implies that one needs to find a globally attracting state $x^*$ where the controllability condition \controllability{x^*} is satisfied to obtain the $\varphi$-irreducible and T-chain property of the Markov kernel.
\begin{theorem}[Practical condition for $\varphi$-irreducibility and T-chain property]
\label{t:practical-condition-ireducibility}
Consider the Markov kernel $P$ defined via \eqref{eq:def-MC}, and suppose \Cref{ass:alpha}-\Cref{ass:sub_differe_S}. Then $P$ is a $\varphi$-irreducible T-chain, and thus every compact set of $\cX$ is petite.
\end{theorem}
\begin{proof}
By \Cref{p:rank-condition-on-all-states}, for any $x\in\cX$, \controllability{x} holds. Then, by \Cref{c:Tchain}, $P$ is a T-chain and by \Cref{t:sufficient-condition-ireducibility}, $P¨$ is $\varphi$-irreducible, and by \cite[Theorem 6.2.5]{meyn2012markov} all compact sets of $\cX$ are petite.
\end{proof}
This latter theorem constitutes the first part of our main result stated in \Cref{theo:1} while the second part relates to the aperiodicity of the kernel which is developed in the next section.

\subsubsection{Aperiodicity}
In this section, we provide conditions for $P$ to be aperiodic. We start with the following characterization which is the counterpart to \Cref{t:sufficient-condition-ireducibility} for a kernel to be  $\varphi$-irreducible aperiodic.

\begin{theorem}
\label{t:sufficient-condition-aperiodicity}
Consider the Markov kernel $P$ defined via \eqref{eq:def-MC}, and suppose \Cref{ass:alpha} and \Cref{ass:F}. If for every $x\in\cX$, \controllability{x} holds, then $P$ is a $\varphi$-irreducible aperiodic Markov kernel if and only if there exists a steadily attracting state.
\end{theorem}

\begin{proof}
First suppose that $P$ is $\varphi$-irreducible and aperiodic. By \Cref{t:sufficient-condition-ireducibility}, there exists a globally attracting state $x^*\in\cX$. Besides, by \Cref{p:implicit-Skx}, there exists an attainable state $y^*$, to which we apply \Cref{p:cycles} (iii), so that there exists a $d$-cycle. However, $P$ is aperiodic, so $d=1$. Thus, by \Cref{p:cycles} (ii), $y^*$ is steadily attracting.

Conversely, suppose that there exists a steadily attracting state $x^*$. By \Cref{p:sufficient-condition-steadily-attracting} (i), $x^*$ is globally attracting, so that by \Cref{t:sufficient-condition-ireducibility}, $P$ is $\varphi$-irreducible. It remains to prove that it is aperiodic. By \Cref{p:steadily-att-state-differentiable}, \rankcondition{x^*} holds for some $k>0$ and $w_{1:k}^*\in\cO_{x^*}^k$. Therefore, we apply \Cref{p:implicit-Skx}, hence $y^*\coloneqq S_{x^*}^k(w^*_{1:k})$ is attainable. Let $U$ be a neighborhood of $x^*$ which satisfies \Cref{p:implicit-Skx} (i). Since $x^*$ is steadily attracting, there exists $T>0$, such that for every $t\geqslant T$, there exists $u_{1:t}\in\cO_{y^*}^t$ such that $z\coloneqq S_{y^*}^t(u_{1:t})\in U$. As $U$ satisfies \Cref{p:implicit-Skx} (i), then there exists $w_{1:k}\in\cO_{z}^k$ such that $S_z^k(w_{1:k})=S_{x^*}^k(w_{1:k}^*)=y^*$. All in all, we have that for every $t\geqslant T$, there exists $w_{1:k+t}\in\cO_{y^*}^{k+t}$ such that $y^*=S_{y^*}^{k+t}(w_{1:k+t})$, hence $y^*\in A_+^{k+t}(y^*)$. By \Cref{p:cycles} (iii), there exists a $1$-cycle, i.e., $P$ is aperiodic.
\end{proof}

We now state our main practical condition to ensure that $P$ is aperiodic.

\begin{theorem}[Practical condition for $\varphi$-irreducibility and aperiodicity]
\label{t:practical-condition-aperiodicity}
Consider the Markov kernel $P$ defined via \eqref{eq:def-MC}, and suppose \Cref{ass:alpha}-\Cref{ass:F} and \Cref{ass:sub_differe_S_steadily}. Then $P$ is a $\varphi$-irreducible aperiodic T-chain, and every compact set of $\cX$ is small.
\end{theorem}

\begin{proof}
By \Cref{p:rank-condition-on-all-states}, \controllability{x} holds for any $x\in\cX$. Thus, by \Cref{t:practical-condition-ireducibility,t:sufficient-condition-aperiodicity}, we have that $P$ is a $\varphi$-irreducible aperiodic T-chain for which compact sets of $\cX$ are petite. Note that, by \cite[Theorem~5.5.7]{meyn2012markov}, any petite set is small.
\end{proof}


\subsection{Proofs for the application to CMA-ES}\label{proof-for-CMA}

\begin{proof}[Proof of \Cref{p:steadily-attracting-cma} (i)]
By \Cref{c:sufficient-condition-steadily-attracting}, it is sufficient to find, for any $\theta_0=(z_0,\ncovmat_0)\in\cX$, a sequence $\seq{v_k}{k\geqslant1}$ such that $v_{1:k}\in\overline{\cO_{\theta_0}^k}$ for every $k\geqslant1$, and
$
\lim_{k\to\infty} S_{\theta_0}^k(v_{1:k}) = (0,I_d) .
$
Since $f$ has Lebesgue negligible level sets, we have, by \Cref{p:density-cma}, that for every $(z,\ncovmat)\in\cX$, and for every $u\in\bR^d$, the element $v=(u,\dots,u)$ of $\cW=\bR^{d\mu}$ belongs to $\overline{\cO_{(z,\ncovmat)}^1}$. 

Then, set 
$
v_1^1=\dots=v^\mu_1=-\ncovmat_0^{-1/2}z_0
$
so that $v_1=(v_1^1,\dots,v^\mu_1)\in\overline{\cO_{\theta_0}^1}$ and $z_1=0$.
Note that $S^1(v_1)=\theta_1\coloneqq (0,\ncovmat_1)$, for some $\ncovmat_1\in\Sdpp$.
Next, consider $(e_1,\dots,e_d)$ an orthogonal basis of eigenvectors of the positive definite matrix $\ncovmat_1$, with $\ncovmat_1e_i=\lambda_i(\ncovmat_1)e_i$, where $\lambda_i(\cdot)$ denotes the function that maps a symmetric matrix to its $i$-th largest eigenvalue (counted with multiplicity). 

Let then $\kappa\geqslant0$ and define
$
v^1_2=\dots=v^\mu_2 = -\ncovmat_1^{-1/2}\kappa e_2
$
and $\theta_2=(z_2,\ncovmat_2)=S^1_{\theta_1}(v_2)$. Then let
$
v^1_3=\dots=v^\mu_3 = -\ncovmat_2^{-1/2}z_2 
$
and $\theta_3=(z_3,\ncovmat_3)=S^1_{\theta_2}(v_3)$.
Then, we have $z_2 = \kappa e_2$ and $z_3=0$. Besides,
$
\ncovmat_2 = r_2^{-1}\times \left( (1-c) \ncovmat_1 + c \kappa^2 e_2e_2^\top \right)
$
with $r_2=\det((1-c) \ncovmat_1 + c \kappa^2 e_2e_2^\top)^{1/d}$ depends continuously on the choice of $\kappa\geqslant0$. Moreover, we have $1-c\leqslant r_2\leqslant 1-c+c\kappa^2$. Then,
\begin{equation}
	r_3	\ncovmat_3 =  \ncovmat_1 + c(1-c)^{-2}\kappa^2\times(r_2+1-c) e_2e_2^\top \eqqcolon K_3 
\end{equation}
for some $r_3>0$.
However, the eigenvalue of the matrix $K_3$ associated to the eigenvector $e_1$ equals $\lambda_1(\ncovmat_1)$ for any value of $\kappa$, while the eigenvalue of $K_3$ associated to the eigenvector $e_2$ depends continuously on $\kappa\geqslant0$ and tends to $+\infty$ when $\kappa\to\infty$ and to $\lambda_2(\ncovmat_1)\leqslant\lambda_1(\ncovmat_1)$ when $\kappa\to0$. Hence, there exists a value of $\kappa\geqslant0$ such that the eigenvalues of $K_3$ associated respectively to the eigenvectors $e_1$ and $e_2$ are equal. Setting $\kappa$ to this value, we get then that $\lambda_1(\ncovmat_3)=\lambda_2(\ncovmat_3)$.

Repeating eventually these steps $(d-1)$ times, we find $v_4,\dots,v_{1+2(d-1)}\in\cW$ such that, denoting $\theta_k=(z_k,\ncovmat_k)=S^k_{\theta_0}(v_{1:k})$, $z_{1+2(d-1)}=0$, and $\lambda_1(\ncovmat_{1+2(d-1)})=\dots=\lambda_d(\ncovmat_{1+2(d-1)})$, with $v_k\in\overline{\cO_{\theta_{k-1}}^1}$ for each $k>0$. However, $\det(\ncovmat_k)=1$ for every $k\in\bN$, thus $\ncovmat_{1+2(d-1)}=I_d$.

For the next steps $k\geqslant 2d$, we choose $v_k=0\in\overline{\cO^{1}_{\theta_{k-1}}}$, so that, by induction, we obtain $\theta_k=(0,I_d)$. By \Cref{c:sufficient-condition-steadily-attracting}, we find that $(0,I_d)$ is a steadily attracting state.
\end{proof}

\begin{proof}[Proof of \Cref{p:steadily-attracting-cma} (ii)]
Let $k>0$ and $v_{1:k}\in\overline{\cO^k_{\theta_0}}$, with $\theta_0=(z_0,\ncovmat_0)\coloneqq (0,I_d)$.
We find here values for $k>0$ and $v_1,\dots,v_{k}\in\cW$ such that the map
$$
\mathcal{D}S^{k}_{(0,I_d)}\left(v_{1:k}\right) \colon T_{(v_{1:k})}\cW^k \to T_{S^k_{(0,I_d)}(v_{1:k})}\cX
$$
is full-rank, i.e., is surjective. We remind that $\cW=(\bR^d)^\mu$, hence
$
T_{(v_{1:k})}\cW^k =\cW^k = \bR^{d\times\mu\times k} .
$
Moreover, we have $\cX= \bR^d\times \det^{-1}(\{1\})$, therefore, by \cite[Proposition 5.38]{lee2012introduction}
$$
T_{S^k_{(0,I_d)}(v_{1:k})}\cX = \bR^d \times \ker ~ \cD\det\left(\ncovmat_k\right) ,
$$
where $\theta_t=(z_t,\ncovmat_t)= S_{\theta_0}^t(v_{1:t})$ for each $t=0,\dots,k$.

We define then inductively the covariance matrix before normalization as
$$
K_{t+1} = (1-c) K_t + c \sqrt{K_t} \sum_{i=1}^\mu w_i \left(v_{t+1}^i\right)\left(v_{t+1}^i\right)^\top \sqrt{K_t}
$$
with $K_0=\ncovmat_0=I_d$, so that, by induction, we have for every $t=0,\dots,k$ ,
\begin{math}
	\ncovmat_t = \frac{K_t}{\det(K_t)^{1/d}} .
\end{math}
Let us introduce (small) perturbations $h_t=(h^1_t,\dots,h^\mu_t)\in\cW$ for $t=1,\dots,k$, and let us denote the perturbed process as
$$
\theta_t^h=(z_t^h, \ncovmat_t^h) = S^t_{\theta_0}(v_{1:t}+h_{1:t}) .
$$
Define  $K_t^h\in\Sdpp$ similarly.
Set $k_0=d(d+1)/2$ the dimension of $\Sd$, and set $k=k_0(k_0-1)+1$. Then, set $v_1,\dots,v_{k_0}$ as follows. Define $\psi_1,\dots,\psi_{k_0}$ nonzero vectors of $\bR^d$, such that $(\psi_1\psi_1^\top,\dots,\psi_{k_0}\psi_{k_0}^\top)$ forms a basis of $\Sd$.

For $t=1,\dots,k_0$, using \Cref{p:density-cma}, we set
$
v_{t} = (K_{t-1}^{-1/2}\psi_t,\dots,K_{t-1}^{-1/2}\psi_t)\in \overline{\cO^1_{\theta_{t-1}}} ,
$
so that
$
K_{t} = (1-c) K_{t-1} + c \psi_t\psi_t^\top .
$
\newcommand{\h}{\varepsilon}
Fix then $\kappa_{t}^1\in\bR$, and let $\h_1>0$ be an arbitrary small positive quantity. Set
$
h_t^1=\dots=h_t^\mu = \frac{1}{2}\kappa_t^1 \h_1 K_{t-1}^{-1/2}\psi_t
$
and then
$$
K_{t}^h =(1-c) K_{t-1}^h + c \psi_t\psi_t^\top +  \h_1 \kappa_t^1 c \psi_t\psi_t^\top + \h_1 A_t^1(\h_1)
,
$$
where $A_t^1(\h_1)\in\Sd$ tends to $0$ when $\varepsilon_1\to0$.
Then, we get by induction,
$$
K_{k_0}^h = K_{k_0} + \h_1 \sum_{t=1}^{k_0} \kappa_t^1 (1-c)^{k_0-t}c \psi_t\psi_t^\top +\h_1 A_{k_0}^1(\h_1) .
$$
Likewise, $A_{k_0}^1(\h_1)$ defines a symmetric matrix which then tends to $0$ when $\h_1$ tends to $0$.
Repeat these steps $k_0-1$ times with $\h_2,\dots,\h_{k_0-1}>0$ instead of $\h_1>0$ and $\kappa_t^2,\dots,\kappa_t^{k_0-1}\in\bR$ instead of $\kappa_t^1\in\bR$. All in all, we have finally, since $k=k_0(k_0-1)+1$,
$$
K_{k-1}^h=K_{k-1} + \sum_{s=1}^{k_0-1} \left[ \h_s  \sum_{t=1}^{k_0} \kappa_t^s (1-c)^{k_0(k_0-1)-sk_0+k_0-t}c \psi_t\psi_t^\top + \h_s A_{k-1}^s(\h_s) \right] .
$$
Again, for each $s=1,\dots,k_0-1$, $A_{k-1}^s(\h_s)$ defines a symmetric matrix which tends to $0$ when $\h_s$ tends to $0$.
Now, consider $(S_1,\dots,S_{k_0-1})$ a basis of $\ker\mathcal{D}\det(\ncovmat_{k-1})$. 
For $s=1,\dots,k_0-1$, we set now the real values $\kappa_t^s$, $t=1,\dots,k_0$ such that we have
$$
\sum_{t=1}^{k_0} \kappa_t^s (1-c)^{k-1-sk_0+k_0-t}c \psi_s\psi_s^\top = S_s .
$$
This is possible since $(\psi_1\psi_1^\top,\dots,\psi_{k_0}\psi_{k_0}^\top)$ is a basis of $\Sd$.
Then,
$$
K_{k-1}^h=K_{k-1} + \sum_{t=1}^{k_0-1} \h_t S_t + \h_t A_{k-1}^t(\h_t) .
$$
Yet, since $S_t\in\ker\cD\det(\ncovmat_{k-1})$, we have then
\begin{align*}
	\ncovmat_{k-1}^h  = \frac{K_{k-1}^h}{\det(K_{k-1}^h)^{1/d}} 
	& = \frac{K_{k-1} + \sum_{t=1}^{k_0-1} \h_t S_t + \h_t A_{k-1}^t(\h_t)}{\det\left(K_{k-1} + \sum_{t=1}^{k_0-1} \h_t S_t + \h_t A_{k-1}^t(\h_t)\right)^{1/d}} \\
	& = \ncovmat_{k-1} + r\sum_{t=1}^{k_0-1} \h_t S_t + \h_t B^t(\h_t) ,
\end{align*}
where we set $r=\det(K_{k-1})^{-1/d}$, and the symmetric matrices $B^t(\h_t)$ tend to $0$ when $\h_t\to0$, for $t=1,\dots,k_0-1$.
Lastly, set $v_{k}=0\in\overline{\cO^1_{\theta_{k-1}}}$, and let $h_{k}^1=\dots=h_{k}^\mu = \ncovmat_{k-1}^{-1/2} \xi_{k_0}$, for some arbitrary small vector $\xi_{k_0}\in\bR^d$. Then, $\ncovmat_k=\ncovmat_{k-1}$ and
$$
z_{k}^h = z_{k} + (1-c)^{-1/2}\xi_{k_0} + l(\h_1,\dots,\h_{k_0-1}) +\|(\h_1,\dots,\h_{k_0-1},\xi_{k_0})\| h(\h_1,\dots,\h_{k_0-1},\xi_{k_0}) ,
$$
where the map $l\colon \bR^{k_0-1}\to\bR^d$ is linear and the quantity $h(\h_1,\dots,\h_{k_0-1},\xi_{k_0})$ tends to $0$ when $\|(\h_1,\dots,\h_{k_0-1},\xi_{k_0})\|$ to $0$.
Furthermore,
$$
\ncovmat_k^h = \ncovmat_k +  r\sum_{t=1}^{k_0-1} \h_t S_t + \h_t B^t(\h_t) + c \xi_{k_0}\xi_{k_0}^\top .
$$
Finally,
\begin{equation}
	\cfrac{S_{(0,I_d)}^{k}(v_{1:k}+h_{1:k}) - S_{(0,I_d)}^{k}(v_{1:k}) - \left( \begin{array}{c}
			(1-c)^{-1/2}\xi_{k_0} + l(\h_1,\dots,\h_{k_0-1}) \\
			r\sum_{t=1}^{k_0-1} \h_t S_t
		\end{array} \right) }{\|(\h_1,\dots,\h_{k_0-1},\xi_{k_0})\|} 
\end{equation}
tends to $0$ when $\|(\h_1,\dots,\h_{k_0-1},\xi_{k_0})\|\to0$.
Therefore,
\begin{equation}
	\cD S_{(0,I_d)}^{k}(v_{1:k})h_{1:k} = \left( \begin{array}{c}
		(1-c)^{-1/2}\xi_{k_0} + l(\h_1,\dots,\h_{k_0-1}) \\
		r\sum_{t=1}^{k_0-1} \h_t S_t
	\end{array} \right)
\end{equation}
defines a surjective map from $\cW^{k}$ to $\bR^d\times\ker\cD \det (\ncovmat_k)$.
Indeed, if $H_{\ncovmat} \in \ker\cD \det (\ncovmat_k)$ and $h_z\in\bR^d$, then there exist $\h_1,\dots,\h_{k_0-1}\in\bR$ such that 
$r\sum_{t=1}^{k_0-1} \h_t S_t = H_{\ncovmat}$,
and then there exists $\xi_{k_0}\in\bR^d$ such that 
$\cD S_{(0,I_d)}^{k}(v_{1:k})h_{1:k} = (h_z;H_{\ncovmat})$.
\end{proof}

\bibliographystyle{plain}
\bibliography{biblio}

\appendix

\section{Background on manifolds}
We recall below basics of differential geometry needed in the present paper. We refer to \cite{amann2009analysis} for more details.	
\begin{definition}[Manifolds]
	\label{def:manifolds}
	A topological space $\msx$ is said to be a topological manifold of dimension $n$ if it is a second countable Hausdorff space that is locally Euclidean of dimension $n$.
	
	Note that $\cX$ is said to be a Hausdorff space if for every pair of distinct points $x,y\in\cX$, there exist neighborhoods $U$ of $x$ and $V$ of $y$ that are disjoint. Moreover, $\cX$ is said to be second countable if there exists a countable basis, that is, a countable collection $\cB$ of open subsets of $\cX$ such that any open subset of $\cX$ can be written as the union of sets in $\cB$.
	
	Finally, $\cX$ is locally Euclidean when for every $x\in\cX$, there exists a neighborhood $U$ of $x$, an open set $V$ of $\bR^n$ and a homeormophism (i.e., a continuous bijection with a continuous reciprocal function) $\varphi\colon U\to V$. We call $(\varphi,U)$ a chart around $x$.
	
	Besides, a manifold $\cX$ is said to be smooth if it is topological, locally Euclidean, and if \del{for }every charts $(\varphi,U)$ and $(\psi,V)$ around any point $x\in\cX$ are such that $\varphi\circ\psi^{-1}$ is continuously differentiable.
\end{definition}
We call throughout the paper a manifold a smooth manifold.

Given $\cX$ a ($n$-dimensional) manifold, and $x\in\cX$, we denote by $T_x\cX$ the tangent space of $\cX$ in $x$.
We refer to \cite[Chapter XI]{amann2009analysis} or to \cite[Chapter 1, Section 2]{guillemin2010differential} for a formal definition of tangent spaces.

We introduce the measurability on a smooth manifold via the following definition. We refer to \cite[Chapter XII]{amann2009analysis} for further details.
\begin{definition}
	A subset $\msa\subset\cX$ is said to be measurable if for all $x\in\cX$, there exists a chart around $x$ denoted $(\varphi,U)$ such that $\varphi(\msa\cap U)$ is measurable (in $\bR^n$).
\end{definition}

\newcounter{defs}
\setcounter{defs}{\value{definition}}

\section{Clarke's generalized derivative of locally Lipschitz functions on manifolds}\label{sec:appendix-clarke}
\setcounter{definition}{\value{defs}}

Clarke's generalized Jacobian is defined for locally Lipschitz functions $g\colon\bR^n\to\bR^m$~\cite{clarke1990optimization}. We define here the Clarke's derivative for locally Lipschitz functions $f\colon\cX\to\cY$ where $\cX$ and $\cY$ are smooth manifolds.
First, let us define formally what a locally Lipschitz function between manifolds is.	
\begin{definition}
	\label{def:locallyLipschitz}
	Let $\cX$ and $\cY$ be two manifolds, equipped with their distance functions $d_{\cX}$ and $d_{\cY}$ respectively, and $f\colon\cX\to\cY$ a function. 
	\begin{enumerate}
		\item[(i)] $f$ is said to be Lipschitz if there exists $L>0$ such that for all $x,y\in\cX$ we have $d_{\cY}(f(x),f(y)) \leqslant L\times d_{\cX}(x,y)$.
		\item[(ii)] $f$ is said to be locally Lipschitz at $x\in\cX$ if there exists a neighborhood $U$ of $x$ in $\cX$ such that the restriction of $f$ to $U$ is Lipschitz.
	\end{enumerate}
\end{definition}
As stated below, a function is locally Lipschitz if and only if it is locally Lipschitz in the charts.
\begin{proposition}
	\label{p:lischitz-manifold}
	If $f\colon\cX\to\cY$ is locally Lipschitz at $x\in\cX$, then for all local charts $(\varphi,U)$ of $\cX$ around $x$ and $(\psi,V)$ of $\cY$ around $f(x)$, the function $\psi\circ f\circ\varphi^{-1}$ is locally Lipschitz at $\varphi(x)$.
\end{proposition}
\begin{proof}
	See that both $\varphi^{-1}$ and $\psi$ are $\cC^1$ hence are locally Lipschitz at all points of their domains. By composition we find that $\psi\circ f\circ \varphi^{-1}$ is locally Lipschitz at $\varphi(x)$.
\end{proof}
Rademacher's theorem~\cite[Theorem 3.2]{gariepy2015measure}, states that a locally Lipschitz function is almost everywhere differentiable. This is easily extended to locally Lipschitz functions on manifolds.
\begin{corollary}[Rademacher's theorem]
	\label{c:rademacher}
	Let $\zeta_{\cX}$ be a measure on $\cX$, which is locally equivalent to the Lebesgue measure, that is, for any measurable subset $A$ of $\cX$, then $\zeta_{\cX}(A)=0$ if and only if for every charts $(\varphi,U)$ of $\cX$, $\Leb\circ\varphi(A\cap U)=0$.
	Then, any function $f\colon\cX\to\cY$ locally Lipschitz at every $x\in\cX$, is differentiable $\zeta_{\cX}$-almost everywhere.
\end{corollary}
\begin{proof}
	Consider local charts $(\varphi,U)$ of $\cX$ around $x$ and $(\psi,V)$ of $\cY$ around $f(x)$. Let us prove that for $\zeta_{\cX}$-almost every point $y$ of $U$, $f$ is differentiable at $y$. See that by \Cref{p:lischitz-manifold}, $g\coloneqq\psi\circ f\circ\varphi^{-1}$ is locally Lipschitz on $\varphi(U)$. Thus, by \cite[Theorem 3.2]{gariepy2015measure}, we have that $g$ is differentiable $\Leb$-almost everywhere on $\varphi(U)$. Thus, since $\varphi$ and $\psi^{-1}$ are $\cC^1$, and since the measures $\zeta_{\cX}$ and $\Leb\circ\varphi$ are equivalent on $U$, then $f=\psi^{-1}\circ g\circ \varphi$ is differentiable $\zeta_{\cX}$-almost everywhere.
\end{proof}
We give now the definition of Clarke's Jacobian for locally Lipschitz functions on Euclidean spaces.
\begin{definition}[Clarke's generalized Jacobian]
	\label{def:clarke-jacobian}
	Let $f\colon \bR^n\to\bR^m$ be locally Lipschitz at $x_0\in\bR^n$. Define
	\begin{equation}
		\partial f (x_0) = {\mathrm{conv}}\left\{ \lim_{t\to\infty} \cD f(x_t) \mid x_t\to x_0, ~ f \text{ is differentiable in all } x_t \right\}
	\end{equation}
	where $\cD f(x_t)\in\bR^{n\times m}$ is the Jacobian matrix of $f$ at $x_t$ (when defined) and ${\mathrm{conv}}$ denotes the convex hull.
\end{definition}
We generalize now this definition to locally Lipschitz functions on manifolds.
\begin{definitionProposition}[Clarke's generalized Jacobian on manifolds]
	\label{p:clarke-derivative-manifold-chain-rule-supp}
	Let $\cX$ and $\cY$ be two manifolds.	
	Let $f\colon \cX\to\cY$ be locally Lipschitzian at $x_0\in\cX$. Let $(\varphi,U)$ be a local chart of $\cX$ around $x_0$ and $(\psi,V)$ be a local chart of $\cY$ around $f(x_0)$. Define $g=\psi\circ f\circ \varphi^{-1}$. Then $g\colon\varphi(U)\to\psi(V)$ is locally Lipschitz at $\varphi(x_0)$, and we can define
	\begin{equation}\label{eq:def-jac-manifolds-app}
		\partial f (x_0) = \left\{ \cD\psi^{-1}(g\circ\varphi(x_0))\circ h \circ \cD \varphi (x_0) \mid h\in\partial g (\varphi(x_0)) \right\} .
	\end{equation}
\end{definitionProposition}
\begin{proof}
	The maps $\psi$ and $\varphi^{-1}$ are by definition continuously differentiable, hence are locally Lipschitz. Therefore, by composition, $g$ is locally Lipschitz. Furthermore, note that the expression \eqref{eq:def-jac-manifolds-app} does not depend on the choice of the charts. Indeed, let $(\varphi_1,U)$ and $(\varphi_2,U)$ be two charts of $\cX$ at $x_0$ and $(\psi_1,V)$ and $(\psi_2,V)$ be two charts of $\cY$ at $f(x_0)$, such that $g_1 = \psi_1\circ f\circ \varphi_1^{-1}$ and $g_2=\psi_2\circ f\circ\varphi_2^{-1}$ are well defined. Then, note that $g_2 = \psi_2\circ\psi_1^{-1} \circ g_1 \circ \varphi_1\circ \varphi_2^{-1}$. Apply then the chain rule~\cite[Corollary of Theorem 2.6.6]{clarke1990optimization} to $g_2$ and get
	\begin{multline*}
		\partial g_2 (\varphi_2(x_0)) \\  = \left\{  \cD\psi_2 (f(x_0))\cD\psi_1^{-1}(g_1(\varphi_1(x_0))) \circ h\circ \cD\varphi_1(x_0) \cD\varphi_2^{-1}(\varphi_2(x_0)) \mid h \in \partial g_1(\varphi_1(x_0))\right\} .
	\end{multline*}
	Therefore,
	\begin{multline*}
		\left\{ \cD\psi_1^{-1}(g_1\circ\varphi_1(x_0))\circ h \circ \cD \varphi_1 (x_0) \mid h\in\partial g_1 (\varphi_1(x_0)) \right\} \\ = \left\{ \cD\psi_2^{-1}(g_2\circ\varphi_2(x_0))\circ h \circ \cD \varphi_2 (x_0) \mid h\in\partial g_2 (\varphi_2(x_0)) \right\} .
	\end{multline*}
\end{proof}
We also state the next result, which would be useful to prove \Cref{p:lsc-maxrank}.
\begin{proposition}
	\label{p:compact-clarke-gradient}
	If $f\colon \cX\to\cY$ is locally Lipschitzian at $x_0\in\cX$, then $\partial f (x_0)$ is nonempty, compact and convex.
\end{proposition}

\begin{proof}
	This follows from \cite[Proposition 2.6.2(a)]{clarke1990optimization} and \Cref{p:clarke-derivative-manifold-chain-rule}.
\end{proof}

We now transpose the uppercontinuity of Clarke's Jacobians to the context of locally Lispchitz functions between manifolds.
\begin{proposition}
	\label{p:usc-clarke-jacobian}
	Let $f\colon\cX\to\cY$ be locally Lipschitz at $x_0$. Then, $\lim_{x\to x_0} \partial f(x) \subset \partial f(x_0)$.
\end{proposition}
\begin{proof}
	Let $(\varphi,U)$ and $(\psi,V)$ be two local charts respectively of $\cX$ and $\cY$ at $x_0$ and $f(x_0)$, such that $\tilde f=\psi\circ f\circ \varphi^{-1}$ is well defined. Then, by \Cref{p:clarke-derivative-manifold-chain-rule-supp}, we have, for any $x\in U$
	$$
	\partial f (x) = \left\{ \cD\psi^{-1}(\tilde f\circ\varphi(x))\circ h \circ \cD \varphi (x) \mid h\in\partial \tilde f (\varphi(x)) \right\} .
	$$
	By applying \cite[Proposition 2.6.2]{clarke1990optimization} to $\tilde f$, we find that $\lim_{x\to x_0} \partial \tilde f(\varphi(x)) \subset \partial \tilde f(x_0)$, which ends the proof.
\end{proof}

The next proposition is actually a very important requirement for our analysis. It states that if we can find a point for which the generalized differential of a locally Lipschitz function in this point is of maximal rank, then we can find a point closed to it in which the function is differentiable and the derivative is full rank.
\begin{proposition}
	\label{p:lsc-maxrank-supp}
	Suppose that $f\colon \cX\to\cY$ is locally Lipschitzian at $x_0\in\cX$, and that $\partial f (x_0)$ is of maximal rank, i.e., all $h\in\partial f (x_0)$ is of maximal rank. Then, there exists a neighborhood $U$ of $x_0$ such that for all $y\in U$, $\partial f(y)$ is of maximal rank. Moreover, for all neighborhood $V\subset U$ of $x_0$, there exists $y_0\in V$ such that $f$ is differentiable at $y_0$ and $\cD f(y_0)$ is of maximal rank.
\end{proposition}

\begin{proof}
	Let $A=\{h\in\cL(T_{x_0}\cX,T_{f(x_0)}\cY)\mid h \text{ is not of maximal rank}\}$. Since the application $\rank$ is l.s.c., then $A$ is a closed set. By \Cref{p:compact-clarke-gradient}, $\partial f(x_0)$ is compact, and disjoint from $A$ since it is assumed to be of maximal rank. Thus $d(\partial f(x_0), A)>0$, where $d$ is a metric induced by some norm on the affine space $\cL(T_{x_0}\cX,T_{f(x_0)}\cY)$ of finite dimension. Moreover, there exits $h^*\in\partial f(x_0)$ such that for all $h\in\partial f(x_0)$ we have
	$$
	d(h,A)\geqslant d(h^*,A)=d(\partial f(x_0), A) >0 ,
	$$
	By \cite[Proposition 2.6.2(c)]{clarke1990optimization}, there exists a neighborhood $U$ of $x_0$ such that for all $y\in U$, $d(\partial f(y),A) \geqslant d(h^*,A)/2 >0$, thus $\partial f(y)$ is of maximal rank. The second part follows from \Cref{c:rademacher}.
\end{proof}

Next, we state a chain rule for the generalized Jacobian on manifolds.
\begin{proposition}[Chain rule]
	\label{p:chain-rule-manifolds}
	Let $\cX$, $\cY$ and $\cZ$ be three manifolds. If $f\colon\cX\to\cY$ is locally Lipschitz at $x_0\in\cX$, and if $g\colon\cY\to\cZ$ is differentiable at $f(x_0)$, then we have
	\begin{equation*}
		\partial (g\circ f) (x_0) = \{ \cD g(f(x_0)) h \mid h\in\partial f(x_0) \} .
	\end{equation*}
\end{proposition}

\begin{proof}
	Let $(\varphi,U)$, $(\psi,V)$ and $(\nu,W)$ be local charts respectively of $\cX$, $\cY$ and $\cZ$ around $x_0$, $f(x_0)$ and $g\circ f(x_0)$. Define $\tilde f = \psi\circ f\circ \varphi^{-1}\colon\varphi(U)\to\psi(v)$ and $\tilde g = \nu\circ g\circ \psi^{-1}\colon\psi(V)\to\nu(W)$. Then, by \Cref{p:clarke-derivative-manifold-chain-rule-supp}, we obtain
	\begin{align*}
		\partial (g\circ f) (x_0) & = \left\{ \cD \nu^{-1}(g\circ f(x_0)) \circ H \circ \cD\varphi(x_0) \mid H\in\partial (\tilde g\circ \tilde f)(\varphi(x_0)) \right\} .
	\end{align*}
	Now we apply the chain rule from \cite[Corollary of Theorem 2.6.6]{clarke1990optimization}. Since $\tilde g$ is differentiable at $\tilde f(\varphi(x_0))$ and $\tilde f$ is locally Lipschitz at $\varphi(x_0)$, we have then
	\begin{equation*}
		\partial (\tilde g\circ \tilde f) (\varphi(x_0)) = \left\{ \cD \tilde g(\tilde f(\varphi(x_0))) \circ H \mid H\in\partial \tilde f (\varphi(x_0)) \right\} .
	\end{equation*}
	Then, 
	\begin{align*}
		\partial (g\circ f) (x_0) & = \left\{ \cD \nu^{-1}(g\circ f(x_0)) \circ \cD \tilde g(\tilde f(\varphi(x_0))) \circ H \circ \cD\varphi(x_0) \mid H\in\partial \tilde f (\varphi(x_0)) \right\} \\
		& = \left\{ \cD g(f(x_0)) \circ \cD\psi^{-1}(f(x_0)) \circ H \circ \cD\varphi(x_0) \mid H\in\partial \tilde f (\varphi(x_0)) \right\} \\
		& = \left\{ \cD g(f(x_0)) \circ H \mid H\in\partial f (x_0) \right\} ,
	\end{align*}
	the last line being obtained by applying \Cref{p:clarke-derivative-manifold-chain-rule-supp} to $f$.
\end{proof}

Lastly, the next two theorems are extensions of the inverse function theorem and of the implicit function theorem to our context.	
\begin{theorem}[Inverse function theorem]
	\label{t:inverse-function-theorem}
	Let $\cX$ and $\cY$ be two manifolds of dimension $n$.		
	Let $f\colon \cX\to\cY$ be locally Lipschitzian at $x_0\in\cX$. Suppose that $\partial f(x_0)$ is of maximal rank, i.e., for all $h\in\partial f (x_0)$, we have $\rank h=n$. Then, there exist a neighborhood of $x_0$ in $\cX$, a neighborhood $V$ of $f(x_0)$ in $\cY$ and a Lipschitzian function $g\colon V\to U$ such that
	\begin{enumerate}
		\item[(i)] $g(f(u))=u$ for all $u\in U$ ;
		\item[(ii)] $f(g(v))=v$ for all $v\in V$.
	\end{enumerate}
\end{theorem}

\begin{proof}
	Let $(\varphi,U)$ be a local chart of $\cX$ around $x_0$ and $(\psi,V)$ a local chart of $\cY$ around $f(x_0)$. Define then $\tilde f=\psi\circ f \circ \varphi^{-1}$. Since $\partial f(x_0)$ is of maximal rank, by the chain rule, using \Cref{p:clarke-derivative-manifold-chain-rule}, then $\partial \tilde f (\varphi(x_0))$ is of maximal rank. Then, up to taking $U$ and $V$ smaller, by the Inverse function theorem applied to $\tilde f$ as stated in \cite[Theorem 7.1.1]{clarke1990optimization}, then there exists a Lipschitz function $\tilde g \colon \psi(V)\to\varphi (U)$ such that $\tilde g (\tilde f (\tilde u)) = \tilde u$ for $\tilde u\in\varphi(U)$ and $\tilde f (\tilde g (\tilde v)) = \tilde v$ for $\tilde v\in\varphi(V)$. Define then $g = \varphi^{-1}\circ \tilde g \circ \psi\colon U\to V$ to get
	$$
	g(f(u)) = \varphi^{-1} \circ \tilde g \circ \psi (f (u)) = \varphi^{-1} \circ g(\tilde f (\varphi(u)) )= \varphi^{-1}\circ \varphi (u) = u
	$$
	for all $u\in U$, and
	$$
	f(g(v)) = f\circ\varphi^{-1} \circ \tilde g \circ \psi (v) = \psi^{-1} \circ \tilde f ( \tilde g(\psi (v))) = \psi^{-1}\circ \psi(v) =v
	$$
	for all $v\in V$.
\end{proof}

\begin{theorem}[Implicit function theorem]
	\label{t:implicit-function-theorem}
	Let $\cX$, $\cY$ and $\cZ$ be manifolds of dimensions respectively $m$, $k$ and $k$.
	Let $f\colon \cX\times\cY\to\cZ$ be locally Lipschitzian at $(x_0,y_0)\in\cX\times\cY$. Moreover, assume that the partial generalized differential $\partial_y f(x_0,y_0)$ is of maximal rank. Then there exists a neighborhood $U$ of $x_0$ and a Lipschitz function $g\colon U\to \cY$ such that $g(x_0)=y_0$, and for all $x\in U$,
	\begin{equation}
		f(x,g(x)) = f(x_0,y_0) .
	\end{equation}
\end{theorem}

\begin{proof}
	Define $F(x,y)=(x,f(x,y))$ a function $\cX\times\cY\to\cX\times\cZ$, which is locally Lipschitz at $(x_0,y_0)$. Define $n=m+k$ and note that the dimensions of $\cX\times\cY$ and $\cX\times\cZ$ both equal $n$. Besides, since $\partial_y f(x_0,y_0)$ is of maximal rank, we find that $\partial F(x_0,y_0)$ is of maximal rank. Thus we can apply the inverse function theorem to $F$ and find neighborhoods $U$, $V$, and $W$ respectively of $x_0$ in $\cX$, $y_0$ in $\cY$ and $f(x_0,y_0)$ in $\cZ$, as well as a Lipschitz function $G\colon U\times W\to U\times V$ such that for all $(x,z)\in U\times W$ we have
	$$
	F(G(x,z)) = (x,z).
	$$
	Note that then $G(x,z)= (x,\tilde G(x,z))$ for some $\tilde G(x,z)\in V$, so that $f(x, \tilde G(x,z))=z$. Therefore, define $g(x) = \tilde G(x,f(x_0,y_0))$ to get
	$$
	f(x,g(x)) = f(x_0,y_0) .
	$$
\end{proof}

\section{Additional proofs}
\label{sec:additional-proofs}

\begin{proposition}
	\label{p:cma-control-model}
	Consider an objective function $f\colon\bR^d\to\bR$ which is scaling-invariant w.r.t.\ $x^*=0$.			
	Then, the sequence $\{z_k,\ncovmat_k\}_{k\in\bN}$ defined by \eqref{eq:normMC} is a time-homogeneous Markov chain which follows \eqref{eq:def-MC}, with functions $F$ and $\alpha$ defined by \eqref{eq:FCMA} and \eqref{eq:alphaCMA} respectively.
\end{proposition}

\begin{proof}
	Let $i=1,\dots,\lambda$. Then, we have
	$$
	f\left(z_\t+\sqrt{\ncovmat_\t}U_{\t+1}^{i}\right) = f\left( R(C_\t)^{-1/2} \times \left[ m_\t + \sqrt{C_\t} U_{\t+1}^i \right] \right).
	$$
	Since $f$ is scaling-invariant, this implies that the permutation $s_{\t+1}$ satisfies almost surely that
	$$
	f\left(z_\t+\sqrt{\ncovmat_\t}U_{\t+1}^{s_{k+1}(1)}\right) \leqslant \dots \leqslant f\left(z_\t+\sqrt{\ncovmat_\t}U_{\t+1}^{s_{k+1}(\lambda)}\right) .
	$$
	Let $\t\geqslant1$, and observe that
	\begin{multline*}
		K_{\t+1}\coloneqq \ncovmat_\t + \sqrt{\ncovmat_\t} \sum_{i=1}^\mu w_i \left[U_{\t+1}^{s_{\t+1}(i)}\right]\left[U_{\t+1}^{s_{\t+1}(i)}\right]^\top\sqrt{\ncovmat_\t}
		\\
		=R(\ncovmat_\t)^{-1}\ncovmat_{\t+1}=R(C_\t)^{-1}C_{\t+1} .
	\end{multline*}
	
	Since $R=\det^{1/d}(\cdot)$ is (positively) homogeneous
	$
	R(K_{\t+1}) = R(C_\t)^{-1} R(C_{\t+1})
	$.
	Further\-more, we have that
	\begin{align*}
		z_{\t+1} = & \left( R(C_{\t+1}) \right)^{-1/2}  \times m_{\t+1} \\
		= &  R(C_\t)^{1/2} R(C_{\t+1})^{-1/2} R(C_\t)^{-1/2}  \times \left[ m_\t +\sqrt{C_\t}\sum_{i=1}^\mu w_i U_{\t+1}^{s_{\t+1}(i)}\right] \\
		= &  	R(K_{\t+1})^{-1/2} \times \left[ m_\t + \sqrt{\ncovmat_\t} \sum_{i=1}^\mu \wi U_{\t+1}^{s_{t+1}(i)} \right] . 
	\end{align*}
	Moreover, 
	\begin{align*}
		\ncovmat_{\t+1} & = R(C_{\t+1})^{-1} C_{\t+1} 
		= R(K_{\t+1})^{-1} K_{\t+1}
	\end{align*}
	All in all, we have 
	$$
	(z_{\t+1},\ncovmat_{\t+1}) = F\left( \left(z_\t,\ncovmat_\t\right) , \alpha\left( \left(z_\t,\ncovmat_\t\right) , \left( U_{\t+1}^1,\dots,U_{\t+1}^\lambda \right)  \right) \right) ,
	$$
	ending the proof.
\end{proof}

\begin{proof}[Proof of \Cref{p:density-cma} and \Cref{p:density-csa}]
	Let $\theta=(z,\ncovmat)\in\cX$, consider i.i.d.\ random variables $U^1,\dots,U^\lambda\sim\mathcal{N}(0,I_d)$, and let $U=(U^1,\dots,U^\lambda)$. Then $W=\alpha(\theta,U)$ satisfies a.s.
	$$
	W = \sum_{\sigma\in\mathfrak{S}_\lambda} \1\left\{ f\left(z+\sqrt{\ncovmat}U^{\sigma(1)}\right) < \dots < f\left(z+\sqrt{\ncovmat}U^{\sigma(\lambda)}\right) \right\} \times \left(U^{\sigma(1)},\dots,U^{\sigma(\mu)} \right)  
	$$
	where $\mathfrak{S}_\lambda$ is the set of permutations of $\{1,\dots,\lambda\}$.
	
	Hence, by symmetry,
	\begin{multline*}				
		W = \frac{1}{(\lambda-\mu)!}\sum_{\sigma\in\mathfrak{S}_\lambda} \1\left\{ f\left(z+\sqrt{\ncovmat}U^{\sigma(1)}\right) < \dots < f\left(z+\sqrt{\ncovmat}U^{\sigma(\mu)}\right) \right\} \\ \times \prod_{k=\mu+1}^\lambda \1\left\{ f\left(z+\sqrt{\ncovmat}U^{\sigma(\mu)}\right) < f\left(z+\sqrt{\ncovmat}U^{\sigma(k)}\right) \right\} \times \left(U^{\sigma(1)},\dots,U^{\sigma(\mu)} \right)  .
	\end{multline*}
	Let $\eta\colon\cW\to\bR_+$ be a smooth map with compact support. We obtain
	\begin{align*}
		\bE\left[ \eta(W) \right] & =\frac{1}{(\lambda-\mu)!}\sum_{\sigma\in\mathfrak{S}_\lambda}  \int \1\left\{ f\left(z+\sqrt{\ncovmat}u_{\sigma(1)}\right) < \dots < f\left(z+\sqrt{\ncovmat}u_{\sigma(\mu)}\right) \right\} \\
		& \quad \quad\quad \quad\quad \quad\times \prod_{k=\mu+1}^\lambda \1\left\{ f\left(z+\sqrt{\ncovmat}u_{\sigma(\mu)}\right) < f\left(z+\sqrt{\ncovmat}u_{\sigma(k)}\right) \right\} \\
		&  \quad \quad\quad \quad\quad \quad\times \eta\left(u_{\sigma(1)},\dots,u_{\sigma(\mu)} \right) \gamma^d(u_1)\dots\gamma^d(u_\lambda) \mathrm{d}u_1\dots\mathrm{d}u_\lambda .
	\end{align*}
	However, observe that, for each $k=\mu+1,\dots,\lambda$, we have
	$$
	\int \1\left\{ f\left(z+\sqrt{\ncovmat}u_{\sigma(\mu)}\right) < f\left(z+\sqrt{\ncovmat}u_{\sigma(k)}\right) \right\} \gamma^d(u_{\sigma(k)}) \mathrm{d}u_{\sigma(k)} = 1- Q_\theta^f\left(u_{\sigma(\mu)}\right) .
	$$
	We deduce then the desired result. Note that \Cref{p:density-csa} is obtained by taking $\ncovmat=I_d$.
\end{proof}

\begin{proof}[Proof of \Cref{p:characterization-globally-attracting}]
	First observe that \eqref{eq:def-globally-attracting-state} is equivalent to (iii). 
	Indeed, if (iii) holds, then for every $y\in\cX$ there exists a sequence $\{y_k\}_{k\in\bN}$ with $y_k\in A_+^k(y)$, and with a subsequence converging to $x^*$. 
	In that case, for every $T\geqslant1$, and for every neighborhood $N$ of $x^*$, there exist infinitely many indices $k\geqslant T$ such that $y_k \in N$. 
	Thus, every neighborhood $N$ of $x^*$ intersects $\cap_{T\geqslant1}\cup_{k\geqslant T} A_+^k(y)$, which proves that \eqref{eq:def-globally-attracting-state} holds.
	Conversely, assume that \eqref{eq:def-globally-attracting-state} holds.
	Then, for every $T\geqslant1$, there exists $k\geqslant T$ such that $x^*\in \overline{A_+^k(y)}$. Then, consider $y_k\in A_+^k(y)$ such that $\distX(x^*,y_k)\leqslant 1/k$. This proves (iii).

	Now, suppose (iii) and let us prove (i). Let $y\in\cX$. By (iii), we know that there exists a sequence $\{y_k\}_{k>0}$ such that $y_k\in A_+^k(y)$ and with a subsequence which converge to $x^*$. However, $A_+^k(y)\subset A_+(y)$, therefore $\{y_k\}_{k>0}$ is a sequence with values in $A_+(y)$ admitting $x^*$ as an accumulation point, which proves (i).
	
	Next, assume that (i) holds, and let us prove that this implies (ii). Let $y\in\cX$. By (i), $x^*\in\overline{A_+(y)}$. In other words, for any open $U$ of $\cX$ containing $x^*$, we have $U\cap A_+(y)\neq \emptyset$. Let $U$ be such an open subset. Since $A_+(y)=\cup_{k\in\bN}A_+^k(y)$, then there exists $k\in\bN$ such that $A_+^k(y)$ intersects $U$. If $k\geqslant1$, this proves (ii). Else, if $k=0$, we do the same reasoning with $z=S_y^1(w_1)$ for some $w_1\in\cO^1_y$, which proves (ii).
	
	Last, let us prove that (ii) implies (iii). Suppose (ii), let $y\in\cX$, and let $\{k_t\}_{t\geqslant1}$ be an increasing sequence of $\bN_{>0}$ which satisfies that, for every $t\geqslant1$, there exists a $k_t$-steps path from $y$ to $B(x^*,1/t)$. Hence, let $\{y_k\}_{k>0}$ be a sequence such that $y_k\in A_+^k(y)$ for every $k>0$, and with $y_{k_t}\in B(x^*,1/t)$. Then, the subsequence $\{y_{k_t}\}_{t>0}$ converges to $x^*$, proving (iii).
\end{proof}

\begin{proof}[Proof of \Cref{p:equivalence-ksteps}]
	Let $U$ be an open subset of $\cX$, $x\in\cX$ and $k>0$. First let us assume that $P^k(x,U)>0$. However, note that we have
	\begin{equation*}
		P^k(x,U) = \int_{\cO^k_x} \1\{S^k_x(w_{1:k})\in U\} p^k_x(w_{1:k}) \d\zeta_{\cW}(w_1)\dots \d\zeta_{\cW}(w_k) ,
	\end{equation*}
	which implies that there exists $w_{1:k}\in\cO^k_x$ such that $S^k_x(w_{1:k})\in U$, hence $w_{1:k}$ is a $k$-steps path from $x$ to $U$. Conversely, assume that there exists $w_{1:k}$ a $k$-steps path from $x$ to $U$, i.e., that $S^k_x(w_{1:k})\in U$. Since $F$ is continuous, then $S_x^k$ is continuous as well. Therefore there exists an open subset $V$ of $\cW^k$ such that for all $v_{1:k}\in V$, $S^k_x(v_{1:k})$. Then we obtain
	\begin{equation*}
		P^k(x,U) \geqslant \int_{\cO^k_x\cap V} \1\{S^k_x(w_{1:k})\in U\} p^k_x(w_{1:k})d\zeta_{\cW}(w_1)\dots d\zeta_{\cW}(w_k)  >0 ,
	\end{equation*}
	since $\cO^k_x\cap V$ is open by intersection, and $p_x^k$ is l.s.c.
\end{proof}

\begin{proof}[Proof of \Cref{p:sufficient-condition-steadily-attracting}]
	The statement (i) is a consequence of \Cref{p:characterization-globally-attracting} (ii). 
	
	For (ii), let $y\in\cX$, and for every integer $s\geqslant1$, consider the open subset $U_s=B(x^*,1/s)$ of $\cX$. Then, there exists a nondecreasing sequence $\{T_s\}_{s\geqslant1}$, such that for every $k\geqslant T_s$, there exists a $k$-steps path $w_{1:k}^k\in\cO_{y}^k$ from $y$ to $U_s$.
	For $k\in\bN$, define $y_k = S_y^k(w_{1:k}^k)$, and observe that 
	$y_k\in A_+^k(y)$.
	Moreover, we have $y_k\in U_s$ for every $k\in\{T_s,\dots,T_{s+1}-1\}$. 
	Then, the sequence $\{y_k\}_{k\in\bN}$ converges to $x^*$.
	Conversely, suppose that for every $y\in\cX$, there exists a sequence $\{y_k\}_{k\in\bN}$ converging to $x^*$ with $y_k\in\overline{A_+^k(y)}$.
	Hence, for every $k>0$ there exists $z_k \in A_+^k(y) \cap B(y_k,1/k)$. 
	By definition of $A_+^k(y)$, there exists then $w_{1:k}^k\in\cO_{y}^k$ such that $z_k = S_{y}^k(w_{1:k}^k)$.
	Besides, since $y_k$ tends to $x^*$, then $z_k$ tends to $x^*$ as well. Let $U$ be a neighborhood of $x^*$, so that there exists $T\in\bN$ with $z_k\in U$ when $k\geqslant T$.
	Then, for $k\geqslant T$, $w_{1:k}^k$ is a $k$-steps path from $y$ to $U$. Thus $x^*$ is steadily attracting%
	, proving (ii).
	
	For (iii), suppose that there exist $x^*$ a steadily attracting state and $y^*$ a globally attracting state. Let us prove that $y^*$ is steadily attracting. Let $U$ be a neighborhood of $y^*$ in $\cX$, and let $z\in\cX$. Since $y^*$ is globally attracting, there exist $k>0$ and a $k$-steps path $w_{1:k}\in\cO_{x^*}^k$ from $x^*$ to $U$, i.e., such that $S_{x^*}^k(w_{1:k})\in U$. Since $F$ is continuous, then $S_{x^*}^k$ is continuous, and thus there exists a neighborhood $V$ of $x^*$ such that for every $x\in V$, we have $S_x^k(w_{1:k})\in V$. Moreover, $x\mapsto p_x^k(w_{1:k})$ is lower semicontinuous, so up to taking $V$ a smaller neighborhood of $x^*$, we can assume that $w_{1:k}\in\cO_{x}^k$. Last, $x^*$ is steadily attracting, so there exists $T>0$ such that for every $t\geqslant T$, there exists a $t$-steps path $v_{1:t}$ from $z$ to $V$, i.e., $S^t_{z}(v_{1:t}) \in V$. All in all, for every $t\geqslant T$, there exists a $(t+k)$-steps path $[v_{1:t},w_{1:k}]$ from $z$ to $U$, ending the proof.
\end{proof}

\begin{proof}[Proof of \Cref{c:sufficient-condition-steadily-attracting}]
	The inclusion $\{S_x^k(w_{1:k})\mid w_{1:k}\in\overline{\cO_x^k}\}\subset\overline{A_+^k(y)}$ follows directly from the definition of $A_+^k(y)$ and the continuity of $F$. Let $x^*\in\cX$ and assume that (i) $x^*$ is steadily attracting. Then, by definition, for every $x\in\cX$ and every neighborhood $U$ of $x^*$, there exists $T>0$ such that for every $k\geqslant T$ there is a $k$-steps path from $x$ to $U$, hence (iii) holds. 
	
	Next, assume that (iii) for every $x\in\cX$ and every neighborhood $U$ of $x^*$, there exists $T>0$ such that for every $k\geqslant T$ we can find $w_{1:k}\in\overline{\cO_x^k}$ with $S_x^k(w_{1:k})\in U$. Then, as in the previous proof, for every integer $s\geqslant1$, consider the open subset $U_s=B(x^*,1/s)$ of $\cX$. Therefore there exists a nondecreasing sequence $\{T_s\}_{s\geqslant1}$, such that for every $k\geqslant T_s$, there exists $w_{1:k}\in\overline{\cO_x^k}$ with $S_x^k(w_{1:k})\in U_s$.
	So we find a sequence $\{y_k\}_{k\in\bN}$ such that $y_k\in \{ S_x^k (w_{1:k}) \mid w_{1:k} \in\overline{\cO_x^k} \}$ for $k\in\bN$, and with $y_k\in U_s$ for $k\in\{T_s,\dots,T_{s+1}-1\}$, which proves (ii).
	
	Last, observe that the implication `(ii) implies (i)' follows directly from \Cref{p:sufficient-condition-steadily-attracting}(ii) and the inclusion $\{S_x^k(w_{1:k})\mid w_{1:k}\in\overline{\cO_x^k}\}\subset\overline{A_+^k(y)}$.
\end{proof}

\begin{proof}[Proof of \Cref{p:cycles}]
	First, we prove (i). Observe that $E$ is nonempty. Indeed, $x^*$ is attainable, so there exist $a\in\bN^*$ and $w_{1:a}\in\cO_{x^*}^a$ with $x^*=S_{x^*}^a(w_{1:a})$, and so $x^*=S^{ak}_{x^*}(w_{1:a},\dots,w_{1:a})$. Consider now $a$ and $b$ two elements of $E$ and let us prove that $d\coloneqq\gcd(a,b)\in E$. By definition, there exist $T_a,T_b>0$, such that for every $k\geqslant T_a$, $x^*\in A_+^{ak}(x^*)$, and for every $k\geqslant T_b$, $x^*\in A_{+}^{bk}(x^*)$. Let $T\in\bN$ be larger than $a/d$, so that for every $k\geqslant 0$, the Euclidean division of $(T+k)d$ by $a$ provides us $q_a$ and $r$ two integers such that $(T+k)d= q_a a+ r$. Besides, by Bézout's theorem, we find that $r= q_b b$ for some $q_b\in\bN$, hence $(T+k)d =q_a a + q_b b$. However, by definition of $T_a$ and $T_b$, we have $x^*\in A_+^{(q_a+T_a)a}(x^*)$ and $x^*\in A_+^{(q_b+T_b)b}(x^*)$. All in all, $x^*\in A_+^{(T+k)d+T_a a + T_b b}$, proving (i) since $d$ divides $a$ and $b$.
	
	To prove (ii), observe that if $\gcd(E)=1$, then, by (i), we have $1\in E$. Then, there exists $T\in\bN$ such that for all $k\geqslant T$, $x^*\in A^k(x^*)$. Let $y\in\cX$. Since $x^*$ is attainable, there exists $t\in\bN$ such that $x^*\in A^t(y)$, so that for all $k\geqslant T+t$, $x^*\in A^k(y)$. Thus, $x^*$ is steadily attracting.
	
	For (iii), define $d=\gcd(E)$, and let $D_i = \cup_{r\geqslant0} A_+^{rd+1}(x^*)$ for $i\in\{0,\dots,d-1\}$. First observe that the $D_i$ are disjoint sets. Indeed, $A_+^i(x^*)$ intersects $A_+^j(x^*)$ for some integers $i,j$, then there exists $y$ in their intersection. As $x^*$ is attainable, there exists $k>0$ such that $x^*\in A_+^{k}(y)$, hence $x^*\in A_+^{r(k+i)}(x^*)$ and $x^*\in A_+^{r(k+j)}(x^*)$ for all $r\geqslant0$. This implies that $d$ divides both $k+i$ and $k+j$ hence $d$ divides $i-j$. This shows that the sets $D_i$, $i\in\{0,\dots,d-1\}$, are disjoint. Moreover, by construction, we have $P(y,D^{i+1})=1$ for all $i\in \bZ/d\bZ$. Finally, observe that the union of the $D_i$, $i\in\{0,\dots,d-1\}$, is equal to $A_+(x^*)$. Since $P$ is $\varphi$-irreducible, and $P^k(x^*,A_+(x^*))=1$ for all $k\in\bN$, then the support of $\varphi$ is included in $A_+(x^*)$. All in all, we have that $\{D_i\}_{0\leqslant i\leqslant d-1}$ is a $d$-cycle.
\end{proof}

\begin{proof}[Proof of \Cref{p:support-irreducibility-measure}]
	First we prove \eqref{eq:supp-psi=globally-attracting-states}. Let $x^*\in\supp~\psi$, and let $U$ be a neighborhood of $x^*$. Then, $\psi(U)>0$, which implies that for every $y\in\cX$, $\sum_{k\geqslant0} P^k(y,U)>0$. This is true for every neighborhood $U$ of $x^*$, hence, by \Cref{p:equivalence-ksteps} and \Cref{p:characterization-globally-attracting}, $x^*$ is globally attracting.
	
	Conversely, let $x^*\in\cX$ be a globally attracting state. Then, by \Cref{p:equivalence-ksteps} and \Cref{p:characterization-globally-attracting}, for every neighborhood $U$ of $x^*$ and for every $y\in\cX$, there exists $k>0$ such that $P^k(y,U)>0$, hence $\psi(U)>0$. This implies that $x^*\in\supp~\psi$. All in all, we obtain \eqref{eq:supp-psi=globally-attracting-states}.
	
	Now, let us prove \eqref{eq:supp-psi=A+}. Consider $x^*\in\cX$ a globally attracting state, and let $y^*\in\supp~\psi$. By \eqref{eq:supp-psi=globally-attracting-states}, $y^*$ is then a globally attracting state. Therefore, by \Cref{p:characterization-globally-attracting}, $y^*\in\overline{A_+(x^*)}$.
	
	Conversely, let $y^*\in\overline{A_+(x^*)}$ and let us prove that $y^*$ is globally attracting. Let $U$ be a neighborhood of $y^*$, so that $U$ intersects $A_+(x^*)$. This implies that there exist $k\geqslant1$ and $w_{1:k}\in\cO^k_{x^*}$ such that $S_{x^*}^k(w_{1:k})\in U$. Since $F$ is continuous, then $S_x^k(w_{1:k})\in U$ for every $x$ in a neighborhood $V$ of $x^*$. Besides, $x\mapsto p_x^k(w_{1:k})$ is l.s.c., so, up to taking $V$ smaller, we can assume that $w_{1:k}\in\cO_x^k$ for every $x\in V$. Furthermore, $x^*$ is globally attracting, so, for every $z\in\cX$, there exist $t\in\bN$ and $v_{1:t}\in\cO_z^t$ such that $S_z^t(v_{1:t})\in V$, hence such that $[v_{1:t},w_{1:k}]\in\cO_z^{t+k}$ and $S_z^{t+k}(v_{1:t},w_{1:k})$, proving that $y^*$ is globally attracting. This ends the proof, using \eqref{eq:supp-psi=globally-attracting-states}.
\end{proof}

\end{document}

%% file: header.tex
\usepackage[utf8]{inputenc} 
\usepackage[T1]{fontenc}    

\usepackage{wrapfig,lipsum}

\usepackage{paralist}
\usepackage{booktabs}       
\usepackage{nicefrac}       
\usepackage{multirow}
 \usepackage{xargs}

\usepackage{amsmath}
\usepackage{amsthm}
\usepackage{bm}
\allowdisplaybreaks
\usepackage{amssymb,mathrsfs}
\usepackage{amsfonts}
\usepackage{mathtools}
\usepackage{upgreek}
\usepackage{graphicx}
\usepackage{wrapfig}
\usepackage{xcolor}
\usepackage{soul}
\usepackage{pifont}
\usepackage{bbm}
\usepackage{hyperref}
\usepackage{stmaryrd}
\usepackage{array}
\usepackage{enumitem}

\usepackage{tikz}
\usepackage{pgfplots}
\usepackage{aliascnt}
\usepackage[textwidth=30mm]{todonotes}
\usepackage{cleveref}
\usepackage{accents}
\usepackage{etaremune}
\DeclareMathAlphabet{\mathpzc}{OT1}{pzc}{m}{it}

\makeatletter
\theoremstyle{plain}
\newtheorem{theorem}{Theorem}[section]
\crefname{theorem}{theorem}{Theorems}
\Crefname{Theorem}{Theorem}{Theorems}

\newtheorem*{lemma_nonumber*}{Lemma}

\theoremstyle{plain}
\newtheorem{lemma}{Lemma}[section]

\theoremstyle{plain}
\newaliascnt{corollary}{theorem}
\newtheorem{corollary}[corollary]{Corollary}
\aliascntresetthe{corollary}
\crefname{corollary}{corollary}{corollaries}
\Crefname{Corollary}{Corollary}{Corollaries}

\theoremstyle{plain}
\newaliascnt{proposition}{theorem}
\newtheorem{proposition}[proposition]{Proposition}
\aliascntresetthe{proposition}
\crefname{proposition}{proposition}{propositions}
\Crefname{Proposition}{Proposition}{Propositions}

\theoremstyle{remark}
\newaliascnt{definition}{theorem}
\newtheorem{definition}[definition]{Definition}
\aliascntresetthe{definition}
\crefname{definition}{definition}{definitions}
\Crefname{Definition}{Definition}{Definitions}

\newaliascnt{definitionProposition}{theorem}
\newtheorem{definitionProposition}[theorem]{Proposition and Definition}
\aliascntresetthe{definitionProposition}
\crefname{Proposition and Definition}{Proposition and Definition}{Proposition and Definition}
\Crefname{Proposition and Definition}{Proposition and Definition}{Proposition and Definition}

\theoremstyle{remark}
\newaliascnt{remark}{theorem}
\aliascntresetthe{remark}
\crefname{remark}{remark}{remarks}
\Crefname{Remark}{Remark}{Remarks}

\theoremstyle{remark}
\crefname{example}{example}{examples}
\Crefname{Example}{Example}{Examples}

\crefname{figure}{figure}{figures}
\Crefname{Figure}{Figure}{Figures}

\newtheoremstyle{assumption}{\topsep}{\topsep}{\itshape}{0pt}{}{.}{ }{{\bf\thmname{#1}}{\thmnumber{#2}}\textnormal{\thmnote{ (#3)}}}

\theoremstyle{assumption}
\newtheorem{assumptionH}{{H}}

\crefformat{assumptionH}{{\bf H}{\normalfont #2#1#3}}

\numberwithin{equation}{section}

%% file: def_armand.tex
\newcommand{\markupdraft}[2]{
	\ifthenelse{\equal{#1}{display}}{#2}{}
	\ifthenelse{\equal{#1}{color}}{\color{#2}}{}
}

\newcommand{\newcolored}[3][]{{\markupdraft{color}{#2}#3}
	\ifthenelse{\equal{#1}{}}{}{\markupdraft{display}{{\color{yellow!70!black}[#1]}}}} 

\newcommand{\del}[2][]{{\markupdraft{display}{{\color{yellow}[removed: "#2"[#1]]}}}} 

\renewcommand{\del}[2][]{}

\newcommand{\blank}[1]{}

\newcommand{\vertiii}[1]{{\left\vert\kern-0.25ex\left\vert\kern-0.25ex\left\vert #1 
		\right\vert\kern-0.25ex\right\vert\kern-0.25ex\right\vert_2}}

\usepackage{bbm}
\newcommand*{\1}{\mathbbm{1}}

\newcommand{\rank}{\mathrm{rank}\,}

\newcommand{\bR}{\mathbb{R}}

\newcommand{\bN}{\mathbb{N}}
\newcommand{\bE}{\mathbb{E}}

\newcommand{\bZ}{\mathbb{Z}}

\newcommand{\cB}{\mathcal{B}}

\newcommand{\cX}{\mathcal{X}}
\newcommand{\cY}{\mathcal{Y}}
\newcommand{\cZ}{\mathcal{Z}}

\newcommand{\cC}{\mathcal{C}}
\newcommand{\cO}{\mathcal{O}}

\newcommand{\cU}{\mathcal{U}}

\newcommand{\cL}{\mathcal{L}}

\newcommand{\cW}{\mathcal{W}}
\newcommand{\cD}{\mathcal{D}}
\newcommand{\cR}{\mathcal{R}}

\newcommand{\mueff}{\mu_{\mathrm{eff}}}

\usepackage{xparse}

\newcommand{\Sd}{\mathcal{S}^d}

\newcommand{\Sdpp}{\mathcal{S}^d_{++}}

\newcommand{\wi}[1][i]{{w^m_{#1}}}

\NewDocumentCommand{\yt}{ O{t} }{Y_{#1}}
\NewDocumentCommand{\yit}{ O{i} O{t} }{\left[Y_{#2}\right]_{#1}}

\newcommand{\ncovmat}{\boldsymbol{\Sigma}}

\renewcommand{\t}{k} 

\newcommand{\seq}[2]{\{#1\}_{#2}}

\renewcommand{\cX}{\msx}
\renewcommand{\cW}{\msw}
\renewcommand{\cU}{\msu}
\renewcommand{\cY}{\msy}
\renewcommand{\cZ}{\msz}
\renewcommand{\d}{\mathrm{d}}


%% file: def.tex
\def\msw{\mathsf{W}}

\def\msu{\mathsf{U}}

\def\Argmin{\mathrm{Argmin}}

\def\msa{\mathsf{A}}

\def\msu{\mathsf{U}}

\def\msx{\mathsf{X}}
\def\msz{\mathsf{Z}}
\def\msy{\mathsf{Y}}

\def\mcbb{\mathcal{B}}  

\def\mcu{\mathcal{U}}

\def\rset{\mathbb{R}}

\def\nset{\mathbb{N}}

\newcommand{\tvnorm}[1]{\| #1 \|_{\mathrm{TV}}}

\newcommandx{\psr}[3][3=]{\left\langle#1,#2 \right\rangle_{#3}}
\newcommandx{\normr}[2][2=]{ \left\Vert#1 \right\Vert_{#2}}
\newcommandx{\psrLigne}[3][3=]{\langle#1,#2 \rangle_{#3}}
\newcommandx{\normrLigne}[2][2=]{ \Vert#1 \Vert_{#2}}

\newcommandx{\norm}[2][1=]{\ifthenelse{\equal{#1}{}}{\left\Vert #2 \right\Vert}{\left\Vert #2 \right\Vert^{#1}}}

\newcommand\probaMarkovTilde[2][2=]
{\ifthenelse{\equal{#2}{}}{{\widetilde{\mathbb{P}}_{#1}}}{\widetilde{\mathbb{P}}_{#1}\left[ #2\right]}}

\newcommand{\plusinfty}{+\infty}

\def\eqsp{\;}

\newcommand\sequence[3][2=,3=]
{\ifthenelse{\equal{#3}{}}{\ensuremath{\{ #1_{#2}\}}}{\ensuremath{\{ #1_{#2}, \eqsp #2 \in #3 \}}}}
\newcommand\sequenceD[3][2=,3=]
{\ifthenelse{\equal{#3}{}}{\ensuremath{\{ #1_{#2}\}}}{\ensuremath{( #1)_{ #2 \in #3} }}}

\newcommand\sequenceDouble[4][3=,4=]
{\ifthenelse{\equal{#3}{}}{\ensuremath{\{ (#1_{#3},#2_{#3}) \}}}{\ensuremath{\{  (#1_{#3},#2_{#3}), \eqsp #3 \in #4 \}}}}

\def\iid{i.i.d.}

\def\Exp{\mathrm{Exp}}

\newcommand{\beq}{\begin{equation}}
\newcommand{\eeq}{\end{equation}}

\def\Leb{\mathrm{Leb}}

\def\supp{\mathrm{supp}}

\newcommandx{\gperthmc}[2][1=,2=]{\ifthenelse{\equal{#1}{}}{\Xi}{\ifthenelse{\equal{#2}{}}{\Xi_{h,#1}}{\Xi_{#2,#1}}}}
\newcommandx{\Phiverlet}[2][1=,2=]{\ifthenelse{\equal{#1}{}}{\Phi}{\Phi_{#1}^{\circ (#2)}}}
\newcommandx{\gpertub}[2][1=,2=]{\ifthenelse{\equal{#1}{}}{g}{g_{#1}^{#2}}}
\newcommandx{\Phiverletq}[2][1=,2=]{\ifthenelse{\equal{#1}{}}{\widetilde{\Phi}}{\widetilde{\Phi}_{#1}^{\circ (#2)}}}
\newcommandx{\Phiverletqi}[2][1=,2=]{\ifthenelse{\equal{#1}{}}{\bar{\Psi}}{\bar{\Psi}_{#1}^{(#2)}}}
\newcommandx{\Pkerhmc}[2][1=,2=]{\ifthenelse{\equal{#1}{}}{\mathrm{P}}{\mathrm{P}_{#1, #2}}}
\newcommandx{\tPkerhmc}[2][1=,2=]{\ifthenelse{\equal{#1}{}}{\tilde{\mathrm{P}}}{\tilde{\mathrm{P}}_{#1, #2}}}
\newcommandx{\PkerhmcD}[2][1=,2=]{\ifthenelse{\equal{#1}{}}{\mathrm{K}}{\mathrm{K}_{#1, #2}}}

\def\distX{\mathrm{dist}_{\msx}}
\def\distW{\mathrm{dist}_{\msw}}


%% file: auto-reg.tex
	\subsection{Auto-regressive Riemannian functional random walk}\label{sec:random-walk}
		We first analyze a simple example in order to illustrate our results. We chose to have very strong assumptions for the sake of simplicity. We believe however that they can be generalized with more work, in particular that the manifold is Hadamard or that the density $q(\cdot)$ below is positive everywhere.
        
        Consider the process $\{\phi_k\}_{k\in\bN}$ defined on a
        Hadamard manifold $\cX$ by \eqref{eq:generalized-ula}.  As
        already noted, \eqref{eq:generalized-ula} can be rewritten in
        the form \eqref{eq:def-MC} with $F(x,w)= \Exp_x(w)$ and
        $\alpha(x,u)=-\gamma s(x)+ u$. 
Moreover, under appropriate conditions on $s$ and on the distribution of $\{U_{k+1}\}_{k\in\nset}$, we can apply our results:
\begin{theorem}
Assume that $s$ is locally Lipschitz and that $U_{1}$ admits a density $q(\cdot)$
 with respect to the Lebesgue measure which is positive and lower semicontinuous. Then the Markov chain $\{\phi_k\}_{k\in\bN}$ is an irreducible and aperiodic T-chain. In addition, any compact set is small for the corresponding Markov kernel. 
\end{theorem}

Note that this result is not surprising and could be proven directly without relying on our theory. However, it serves as a simple example where we can easily verify the conditions of \Cref{theo:1}.

			\begin{proof}
				We verify that assumptions of \Cref{theo:1} hold.  We first observe
				that $\alpha(x,U_1)$ satisfies \Cref{ass:alpha} since it has a positive lower semicontinuous
				density with respect to the Lebesgue measure.  
				Moreover, \Cref{ass:alpha} holds
				since the density of $\alpha(x,U_1)$ writes as
					$$
				p_x(u) = q(u + \gamma s(x)) \eqsp,
				$$
				and thus $(x,u)\mapsto p_x(u)$ is l.s.c.\ since $q(\cdot)$ and $s(\cdot)$ are l.s.c.\ as well.
				Second, for \Cref{ass:F}, we simply use that the exponential map is smooth~\cite[Proposition 5.7]{lee1997riemannian}.
				Last, for \Cref{ass:sub_differe_S_steadily}, we first prove that every $x^*\in\cX$ is a steadily attracting state.
				Let $x^*\in\cX$.
				For every $x_0\in\cX$, by Hadamard's theorem $\Exp_{x_0}(\cdot)$ is a diffeomorphism, hence is bijective. Therefore, there exists $w_1\in \mathrm T_{x_0}\cX$ such that $x^* = \Exp_{x_0}(w_1)$. Since moreover, $\Exp_{x^*}(0)=x^*$, then, for every $k\geqslant1$, there exists a $k$-steps path $w_{1:k}=(w_1,0,\dots,0)$ between $x_0$ and $x^*$ and thus $x^*$ is steadily attracting. The path $w_{1:k}$ indeed belongs to the control set $\mathcal O_{x_0}^k$ since $p_x(\cdot)$ is positive for every $x$ (since $q(\cdot)$ is positive).
				
				Furthermore, $\Exp_{x^*}(\cdot)$  is a diffeomorphism,  therefore the Jacobian $\mathcal{D}\Exp_{x^*}(0)$ is invertible and thus of maximal rank.
				Therefore the controllability condition \controllability{x^*} holds, which proves \Cref{ass:sub_differe_S_steadily}. 
				The desired result follows then from \Cref{theo:1}.
			\end{proof}


                        

%% file: verifiable_conditions_manifolds.bbl
\begin{thebibliography}{10}

\bibitem{absil2012projection}
P-A Absil and J{\'e}r{\^o}me Malick.
\newblock Projection-like retractions on matrix manifolds.
\newblock {\em SIAM Journal on Optimization}, 22(1):135--158, 2012.

\bibitem{akiba2019optuna}
Takuya Akiba, Shotaro Sano, Toshihiko Yanase, Takeru Ohta, and Masanori Koyama.
\newblock Optuna: {{A Next-generation Hyperparameter Optimization Framework}}.
\newblock In {\em Proceedings of the 25th {{ACM SIGKDD International
  Conference}} on {{Knowledge Discovery}} \& {{Data Mining}}}, {{KDD}} '19,
  pages 2623--2631, New York, NY, USA, July 2019. Association for Computing
  Machinery.

\bibitem{akimoto2010bidirectional}
Youhei Akimoto, Yuichi Nagata, Isao Ono, and Shigenobu Kobayashi.
\newblock Bidirectional {{Relation}} between {{CMA Evolution Strategies}} and
  {{Natural Evolution Strategies}}.
\newblock In {\em Parallel {{Problem Solving}} from {{Nature}}, {{PPSN XI}}},
  Lecture {{Notes}} in {{Computer Science}}, pages 154--163, Berlin,
  Heidelberg, 2010. Springer.

\bibitem{akimoto2012theoretical}
Youhei Akimoto, Yuichi Nagata, Isao Ono, and Shigenobu Kobayashi.
\newblock Theoretical {{Foundation}} for {{CMA-ES}} from {{Information Geometry
  Perspective}}.
\newblock {\em Algorithmica}, 64(4):698--716, December 2012.

\bibitem{amann2009analysis}
Herbert Amann and Joachim Escher.
\newblock {\em Analysis {{III}}}.
\newblock Birkh{\"a}user, Basel, 2009.

\bibitem{an1997note}
H.~Z. An and S.~G. Chen.
\newblock A note on the ergodicity of non-linear autoregressive model.
\newblock {\em Statistics \& Probability Letters}, 34(4):365--372, June 1997.

\bibitem{auger2016these}
Anne Auger.
\newblock Analysis of {{Comparison-based Stochastic Continuous Black-Box
  Optimization Algorithms}}.
\newblock Th{\`e}se d'habilitation {\`a} diriger des recherches, Universit{\'e}
  Paris-Sud, May 2016.

\bibitem{auger2016linear}
Anne Auger and Nikolaus Hansen.
\newblock Linear {{Convergence}} of {{Comparison-based Step-size Adaptive
  Randomized Search}} via {{Stability}} of {{Markov Chains}}.
\newblock {\em SIAM Journal on Optimization}, 26(3):1589--1624, January 2016.

\bibitem{bacak2014convex}
Miroslav Bac{\'a}k.
\newblock {\em Convex analysis and optimization in Hadamard spaces}, volume~22.
\newblock Walter de Gruyter GmbH \& Co KG, 2014.

\bibitem{bharath2025sampling}
Karthik Bharath, Alexander Lewis, Akash Sharma, and Michael~V Tretyakov.
\newblock Sampling and estimation on manifolds using the langevin diffusion.
\newblock {\em Journal of Machine Learning Research}, 26(71):1--50, 2025.

\bibitem{bhattacharya1995geometric}
Rabi Bhattacharya and Chanho Lee.
\newblock On geometric ergodicity of nonlinear autoregressive models.
\newblock {\em Statistics \& Probability Letters}, 22(4):311--315, March 1995.

\bibitem{bieler2014robust}
Jonathan Bieler, Rosamaria Cannavo, Kyle Gustafson, Cedric Gobet, David
  Gatfield, and Felix Naef.
\newblock Robust synchronization of coupled circadian and cell cycle
  oscillators in single mammalian cells.
\newblock {\em Molecular Systems Biology}, 10(7):739, July 2014.

\bibitem{cheng2022efficient}
Xiang Cheng, Jingzhao Zhang, and Suvrit Sra.
\newblock Efficient {{Sampling}} on {{Riemannian Manifolds}} via {{Langevin
  MCMC}}.
\newblock {\em Advances in Neural Information Processing Systems},
  35:5995--6006, December 2022.

\bibitem{chotard2019verifiable}
Alexandre Chotard and Anne Auger.
\newblock Verifiable conditions for the irreducibility and aperiodicity of
  {{Markov}} chains by analyzing underlying deterministic models.
\newblock {\em Bernoulli}, 25(1):112--147, February 2019.

\bibitem{clarke1990optimization}
Frank~H. Clarke.
\newblock {\em Optimization and {{Nonsmooth Analysis}}}.
\newblock SIAM, January 1990.

\bibitem{do1992riemannian}
Manfredo~Perdigao Do~Carmo and J~Flaherty~Francis.
\newblock {\em Riemannian geometry}, volume~2.
\newblock Springer, 1992.

\bibitem{gariepy2015measure}
Lawrence~Craig Evans and Ronald~F Gariepy.
\newblock {\em Measure {{Theory}} and {{Fine Properties}} of {{Functions}},
  {{Revised Edition}}}.
\newblock {Chapman and Hall/CRC}, New York, April 2015.

\bibitem{glynn2017recurrence}
Peter~W. Glynn, Sanatan Rai, and John~E. Glynn.
\newblock {Recurrence classification for a family of non-linear storage
  models}.
\newblock {\em Probability and Mathematical Statistics}, 37(2):337--353, 2017.

\bibitem{guillemin2010differential}
Victor Guillemin and Alan Pollack.
\newblock {\em Differential {{Topology}}}.
\newblock American Mathematical Soc., 2010.

\bibitem{ha2018recurrent}
David Ha and J{\"u}rgen Schmidhuber.
\newblock Recurrent {{World Models Facilitate Policy Evolution}}.
\newblock {\em Advances in neural information processing systems}, 2018.

\bibitem{hajlasz1993change}
Piotr Haj{\l}asz.
\newblock Change of variables formula under minimal assumptions.
\newblock In {\em Colloquium {{Mathematicae}}}, volume~64, pages 93--101, 1993.

\bibitem{hansen2003reducing}
Nikolaus Hansen, Sibylle~D. M{\"u}ller, and Petros Koumoutsakos.
\newblock Reducing the {{Time Complexity}} of the {{Derandomized Evolution
  Strategy}} with {{Covariance Matrix Adaptation}} ({{CMA-ES}}).
\newblock {\em Evolutionary Computation}, 11(1):1--18, March 2003.

\bibitem{hansen2001completely}
Nikolaus Hansen and Andreas Ostermeier.
\newblock Completely {{Derandomized Self-Adaptation}} in {{Evolution
  Strategies}}.
\newblock {\em Evolutionary Computation}, 9(2):159--195, June 2001.

\bibitem{huang2002ode}
Jianyi Huang, Ioannis Kontoyiannis, and Sean~P. Meyn.
\newblock The {{ODE Method}} and {{Spectral Theory}} of {{Markov Operators}}.
\newblock In {\em Stochastic {{Theory}} and {{Control}}}, Lecture {{Notes}} in
  {{Control}} and {{Information Sciences}}, pages 205--221, Berlin, Heidelberg,
  2002. Springer.

\bibitem{ichihara1974classification}
Kanji Ichihara and Hiroshi Kunita.
\newblock A classification of the second order degenerate elliptic operators
  and its probabilistic characterization.
\newblock {\em Zeitschrift f{\"u}r Wahrscheinlichkeitstheorie und Verwandte
  Gebiete}, 30(3):235--254, 1974.

\bibitem{kliemann1987recurrence}
Wolfgang Kliemann.
\newblock Recurrence and {{Invariant Measures}} for {{Degenerate Diffusions}}.
\newblock {\em The Annals of Probability}, 15(2):690--707, 1987.

\bibitem{lee1997riemannian}
John~M. Lee.
\newblock {\em Riemannian {{Manifolds}}}.
\newblock Graduate {{Texts}} in {{Mathematics}}. Springer, New York, NY, 1997.

\bibitem{lee2012introduction}
John~M. Lee.
\newblock {\em Introduction to {{Smooth Manifolds}}}, volume 218 of {\em
  Graduate {{Texts}} in {{Mathematics}}}.
\newblock Springer, New York, NY, 2012.

\bibitem{li2009monotone}
Chong Li, Genaro L{\'o}pez, and Victoria Mart{\'\i}n-M{\'a}rquez.
\newblock Monotone vector fields and the proximal point algorithm on hadamard
  manifolds.
\newblock {\em Journal of the London Mathematical Society}, 79(3):663--683,
  2009.

\bibitem{li2023riemannian}
Mufan~(Bill) Li and Murat~A. Erdogdu.
\newblock Riemannian {{Langevin}} algorithm for solving semidefinite programs.
\newblock {\em Bernoulli}, 29(4):3093--3113, November 2023.

\bibitem{meyn1987new}
S.~Meyn and P.~Caines.
\newblock A new approach to stochastic adaptive control.
\newblock {\em IEEE Transactions on Automatic Control}, 32(3):220--226, March
  1987.

\bibitem{meyn1989stochastic}
S.~P. Meyn and P.~E. Caines.
\newblock Stochastic controllability and stochastic {{Lyapunov}} functions with
  applications to adaptive and nonlinear systems.
\newblock In {\em Stochastic {{Differential Systems}}}, Lecture {{Notes}} in
  {{Control}} and {{Information Sciences}}, pages 235--257, Berlin, Heidelberg,
  1989. Springer.

\bibitem{meyn1991asymptotic}
S.~P. Meyn and P.~E. Caines.
\newblock Asymptotic {{Behavior}} of {{Stochastic Systems Possessing Markovian
  Realizations}}.
\newblock {\em SIAM Journal on Control and Optimization}, 29(3):535--561, May
  1991.

\bibitem{meyn2012markov}
Sean~P. Meyn and Richard~L. Tweedie.
\newblock {\em Markov {{Chains}} and {{Stochastic Stability}}}.
\newblock Springer Science \& Business Media, December 2012.

\bibitem{meyn1993model}
S.P. Meyn and L.J. Brown.
\newblock Model reference adaptive control of time varying and stochastic
  systems.
\newblock {\em IEEE Transactions on Automatic Control}, 38(12):1738--1753,
  December 1993.

\bibitem{mokkadem1987criteres}
Abdelkader Mokkadem.
\newblock {\em {Crit{\`e}res de m{\'e}lange pour des processus stationnaires.
  Estimation sous des hypoth{\`e}ses de m{\'e}lange. Entropie des processus
  lin{\'e}aires}}.
\newblock PhD thesis, Universit{\'e} Paris-Sud, September 1987.

\bibitem{patte2022estimation}
C{\'e}cile Patte, Pierre-Yves Brillet, Catalin Fetita, Jean-Fran{\c c}ois
  Bernaudin, Thomas Gille, Hilario Nunes, Dominique Chapelle, and Martin Genet.
\newblock Estimation of {{Regional Pulmonary Compliance}} in {{Idiopathic
  Pulmonary Fibrosis Based}} on {{Personalized Lung Poromechanical Modeling}}.
\newblock {\em Journal of Biomechanical Engineering}, 144(091008), March 2022.

\bibitem{rechenberg1973evolutionsstrategie}
Ingo Rechenberg.
\newblock {\em {Evolutionsstrategie: Optimierung technischer Systeme nach
  Prinzipien der biologischen Evolution}}.
\newblock Frommann-Holzboog, Stuttgart, Germany, 1973.

\bibitem{rodriguez2006hybrid}
Maria {Rodriguez-Fernandez}, Pedro Mendes, and Julio~R. Banga.
\newblock A hybrid approach for efficient and robust parameter estimation in
  biochemical pathways.
\newblock {\em Biosystems}, 83(2):248--265, February 2006.

\bibitem{stroock1972support}
Daniel~W. Stroock and S.~R.~S. Varadhan.
\newblock On the {{Support}} of {{Diffusion Processes}} with {{Applications}}
  to the {{Strong Maximum Principle}}.
\newblock In {\em Contributions to {{Probability Theory}}}, pages 333--360.
  University of California Press, December 1972.

\bibitem{toure2023global}
Cheikh Toure, Anne Auger, and Nikolaus Hansen.
\newblock Global linear convergence of evolution strategies with recombination
  on scaling-invariant functions.
\newblock {\em Journal of Global Optimization}, 86(1):163--203, May 2023.

\bibitem{toure2021scaling}
Cheikh Toure, Armand Gissler, Anne Auger, and Nikolaus Hansen.
\newblock Scaling-invariant {{Functions}} versus {{Positively Homogeneous
  Functions}}.
\newblock {\em Journal of Optimization Theory and Applications},
  191(1):363--383, October 2021.

\bibitem{yao2000stability}
J.-F. Yao and J.-G. Attali.
\newblock On stability of nonlinear {{AR}} processes with {{Markov}} switching.
\newblock {\em Advances in Applied Probability}, 32(2):394--407, June 2000.

\end{thebibliography}
